\documentclass[3pt]{article}

\usepackage{graphicx}\usepackage{epstopdf}\usepackage{color}
\usepackage{amsmath}\usepackage{amsfonts}\usepackage{amssymb}
 \usepackage{amsthm}

\numberwithin{equation}{section}

\begin{document}

\title{Non-conforming finite element methods for transmission eigenvalue
problem}

\author{ { Yidu Yang, Jiayu Han, Hai Bi} \\\\
{\small School of Mathematics and Computer Science, }\\{\small
Guizhou Normal University,  Guiyang,  $550001$,  China}\\{\small
 ydyang@gznu.edu.cn, hanjiayu126@126.com, bihaimath@gznu.edu.cn}}
\date{~}
\pagestyle{plain} \textwidth 145mm \textheight 215mm \topmargin 0pt
\maketitle

\indent{\bf\small Abstract~:~} {\small The transmission eigenvalue problem is an important and challenging topic arising in the inverse scattering theory.
In this paper, for the Helmholtz transmission eigenvalue problem, we
give a weak formulation which is a nonselfadjoint linear eigenvalue
problem. Based on the weak formulation, we first discuss the
non-conforming finite element approximation, and prove the error
estimates of the discrete eigenvalues obtained by the Adini element,
Morley-Zienkiewicz element, modified-Zienkiewicz element et. al. And
we report some numerical examples 
 to validate the efficiency of our approach
  for solving transmission eigenvalue problem.
}\\
\indent{\bf\small Keywords~:~\scriptsize} {\small
transmission eigenvalue, the weak formulation, non-conforming finite
elements, error estimates.}
\section{Introduction}

\indent The transmission eigenvalue problems have important physical
background, for example, they can be used to obtain estimates for
the material properties of the scattering object
\cite{cakoni1,cakoni2,sun0}. In addition, transmission eigenvalues
have theoretical importance in the uniqueness and reconstruction in
inverse scattering theory \cite{colton1}. Before 2010,
significant progresses of the existence of transmission eigenvalues
and applications have been made (see \cite{cakoni2} and the survey paper
\cite{cakoni3}).\\
\indent In recent years, the computation of transmission eigenvalues has
attracted the attention of many researchers. The first numerical
treatment of the transmission eigenvalue problem appears in
\cite{colton2} where three finite element methods, including the Argyris, continuous and mixed finite element methods,
are proposed for the Helmholtz transmission eigenvalues, and has
been further developed  by
\cite{an1,cakoni5,colton2,ji,ji2,monk,sun1,su4,yang2}. In
particular,
\cite{cakoni5} studied the mixed method using the Argyris conforming
elements and \cite{yang2} the $H^{2}$ conforming finite element
method, and made rigorous error analysis.
 Moreover, based on $H^{2}$ conforming finite
element approximations, the iterative methods in
\cite{sun1} and the multigrid method in \cite{ji2} were proposed for computing
real transmission eigenvalues, and two-grid method in \cite{yang2}
for computing real and complex transmission eigenvalues.
And the
spectral-element method was studied in \cite{an1}. However, to the
best of our knowledge, there has no research on the non-conforming
finite element methods for the transmission
eigenvalues even for   arbitrary nonselfadjoint elliptic  eigenvalue problem.\\
\indent Inspired by the works mentioned above,  we transform the fourth order equation of transmission eigenvalue problem   into a
weak formulation,
which is suitable to nonconforming elements. This formulation is a nonselfadjoint linear eigenvalue
problem (see (\ref{s2.11})) with a selfadjoint, continuous and coercive
sesquilinear form $A(\cdot,\cdot)$. Based on the weak formulation we
build a type of non-conforming finite element discretizations with good algebraic structure,
including the Adini element \cite{adini}, modified-Zienkiewicz
element \cite{wang}, Morley-Zienkiewicz element \cite{shi},
12-parameter triangle plate element, 15-parameter triangle plate
element  et. al. (see \cite{shi}). And we prove the error estimates of the numerical
eigenvalues. The proof difficulty lies in the non-symmetry of
right-hand sides of eigenvalue problem  that involves derivatives. To
overcome this difficulty, based on  Babuska-Osborn spectral approximation theory \cite{babuska}, the new proof method employed in this
paper is to establish a fundamental relationship (\ref{s4.18}) and
use it to prove the optimal error estimates
of non-conforming element eigenvalues. \\
\indent For fourth order equation in $\mathbb{R}^3$, it is difficult to implement conforming elements in $H^2$,
  whereas many non-conforming elements such as the Morley-Zienkiewicz element have had their three dimensional versions at present (e.g., see \cite{shi}). Hence it is an
essential and significant work to study the non-conforming element approximation for transmission  eigenvalues.\\
\indent Our non-conforming finite element discretization is easy to
realize under the package of iFEM \cite{chen} with Matlab. We use
the sparse matrix eigenvalue solver $eigs$ to compute the
numerical eigenvalues, 
and numerical results indicate  that our methods are efficient for
computing real and complex transmission eigenvalues as
expected.\\
\indent In this paper, regarding the basic theory of  finite
element methods, we refer to \cite{babuska,brenner,ciarlet,oden,shi}.\\
\indent Throughout this paper, $C$ denotes a positive constant
independent of $h$, which may not be the same constant in different
places. For simplicity, we use the symbol $a\lesssim b$ to mean that
$a\leq C b$.

\section{The weak formulation and non-conforming element method}

 \indent Consider the Helmholtz transmission eigenvalue
problem: Find $k\in \mathbb{C}$, $w, \sigma\in L^{2}(\Omega)$,
$w-\sigma\in H^{2}(\Omega)$ such that
\begin{eqnarray}\label{s2.1}
&&\Delta w+k^{2}nw=0,~~~in~\Omega,\\\label{s2.2}
 &&\Delta
\sigma+k^{2}\sigma=0,~~~in~ \Omega,\\\label{s2.3}
 &&w-\sigma=0,~~~on~ \partial
\Omega,\\\label{s2.4}
 &&\frac{\partial w}{\partial \nu}-\frac{\partial
\sigma}{\partial \nu}=0,~~~ on~\partial \Omega,
\end{eqnarray}
where $\Omega\subset \mathbb{R}^{d}$ (d=2,3) is a bounded simply
connected inhomogeneous medium,
 $\nu$ is the unit outward normal to $\partial \Omega$ and the index of refraction $n(x)$
is positive.\\
\indent Let $W^{s,p}(\Omega)$ denote the usual Sobolev space with
norm $\|\cdot\|_{s,p}$,
 $H^{s}(\Omega)=W^{s,2}(\Omega)$, and $\|\cdot\|_{s,2}=\|\cdot\|_{s}$,
$H^{0}(\Omega)=L^{2}(\Omega)$ with the inner product
$(u,v)_{0}=\int\limits_{\Omega}u\overline{v}dx$. Denote $
H_{0}^{2}(\Omega)=\{v\in H^{2}(\Omega): v|_{\partial
\Omega}=\frac{\partial v}{\partial \nu}|_{\partial \Omega}=0\}. $
 Let $H^{-s}(\Omega)$ be the
``negative space", with  norm given by
\begin{eqnarray*}
\|v\|_{-s}=\sup\limits_{0\not=f\in
H_0^{s}(\Omega)}\frac{|(v,f)_{0}|}{\|f\|_{s}}.
\end{eqnarray*}
It is clear that for any real functions $v_{1}$ and $v_{2}$, norms
$\|v_{1}+iv_{2}\|_{s,p}$ and $\|v_{1}\|_{s,p}+\|v_{2}\|_{s,p}$ are
equivalent in $W^{s,p}(\Omega)$, and norms $\|v_{1}+iv_{2}\|_{-s}$ and $\|v_{1}\|_{-s}+\|v_{2}\|_{-s}$ are equivalent in $H^{-s}(\Omega)$.\\
\indent Define Hilbert space $\mathbf{H}=H_{0}^{2}(\Omega)\times
L^{2}(\Omega)$ with norm
$\|(v,z)\|_{\mathbf{H}}=\|v\|_{2}+\|z\|_{0}$, and define
$\mathbf{H}_{1}=H_0^{1}(\Omega)\times H^{-1}(\Omega)$ with norm
$\|(v,z)\|_{\mathbf{H}_{1}}=\|v\|_{1}+\|z\|_{-1}$.\\
Since $L^{2}(\Omega)\hookrightarrow H^{-1}(\Omega)$ compactly (see
pp.31-39 in \cite{berezanskii}) and $H^{2}(\Omega)\hookrightarrow
H^{1}(\Omega)$ compactly, $\mathbf{H}\hookrightarrow\mathbf{H}_{1}$ compactly. \\
\indent In this paper, we suppose that $n=n(x)\in
L^{\infty}(\Omega)$ satisfying either one of the following
assumptions
\begin{eqnarray*}
&&(C1)~~1+\delta\leq \inf_{\Omega}n(x)\leq n(x)\leq
\sup_{\Omega}n(x)<\infty,\\
&&(C2)~~0< \inf_{\Omega}n(x)\leq n(x)\leq \sup_{\Omega}n(x)<1-\beta,
\end{eqnarray*}
for some constant $\delta>0$ or $\beta>0$.\\
\indent From \cite{cakoni3,rynne} we know that the problem
 (\ref{s2.1})-(\ref{s2.4}) can be written as an equivalent  fourth order
equation for $u=w-\sigma\in H_{0}^{2}(\Omega)$:
\begin{eqnarray*}
(\Delta+k^{2}n)\frac{1}{n-1}(\Delta+k^{2})u=0,
\end{eqnarray*}
i.e.,
\begin{eqnarray}\label{s2.5}
\Delta(\frac{1}{n-1}\Delta u)=-k^{2}\frac{n}{n-1}\Delta
u-k^{2}\Delta (\frac{1}{n-1}u)-k^{4}\frac{n}{n-1}u.
\end{eqnarray}
Then the weak formulation for the transmission eigenvalue problem
(\ref{s2.1})-(\ref{s2.4}) can be stated as follows: Find $k\in
\mathbb{C}$, $u\in H_{0}^{2}(\Omega)$ such that
\begin{eqnarray}\label{s2.6}
&&(\frac{1}{n-1}\Delta u,\Delta v )_{0}=k^{2}(\nabla u,
\nabla(\frac{n}{n-1}v)_{0}+k^{2}(\nabla (\frac{1}{n-1}u),
\nabla v)_{0}\nonumber\\
&&~~~~~~-k^{4}(\frac{n}{n-1}u, v)_{0},~~~\forall v
\in~H_{0}^{2}(\Omega).
\end{eqnarray}
Introduce an auxiliary variable
\begin{eqnarray}\label{s2.7}
\omega=k^{2}u,
\end{eqnarray}
then
\begin{eqnarray}\label{s2.8}
(\omega, z)_{0}=k^{2}(u, z)_{0},~~~\forall z\in L^{2}(\Omega).
\end{eqnarray}
Thus, combining (\ref{s2.6}) and (\ref{s2.8}), we arrive at a linear
weak formulation: Find $(k^{2}, u, \omega)\in \mathbb{C}\times
H_{0}^{2}(\Omega)\times L^{2}(\Omega)$ such that
\begin{eqnarray}\label{s2.9}
&&(\frac{1}{n-1}\Delta u, \Delta v)_{0}=k^{2}
(\nabla(\frac{1}{n-1}u),
\nabla v)_{0}\nonumber\\
&&~~~~~~+k^{2}(\nabla u, \nabla(\frac{n}{n-1}v))_{0}-k^{2}
(\frac{n}{n-1}\omega, v )_{0},~~~\forall v
\in~H_{0}^{2}(\Omega),\\\label{s2.10}
 &&(\omega,
z)_{0}=k^{2}(u,z)_{0},~~~\forall z\in L^{2}(\Omega).
\end{eqnarray}
With this weak formulation, we have discussed the conforming finite
element approximations (see \cite{yang2}). However, for the non-conforming element
approximations, the weak formulation can not guarantee that the
discrete bilinear form satisfies the uniform $\mathbf{H}_{h}$-ellipticity
(see Remark 49.1 in \cite{ciarlet}). To study the non-conforming
element approximations,
next we will give a new weak formulation referring to the weak formulation of the
plate problem ( see (49.3) in \cite{ciarlet}).\\
\indent  If $(C1)$ holds, let
\begin{eqnarray*}
&&A((u,\omega),(v,z))=((\frac{1}{n-1}-\mu_{1})\Delta u, \Delta v)_{0}+(\mu_{1}\Delta u, \Delta v)_{0}+(\omega, z)_{0}\\
&&=((\frac{1}{n-1}-\mu_{1})\Delta u, \Delta
v)_{0}+\mu_{1}\int\limits_{\Omega}\sum\limits_{1\leq i,j\leq d}
\frac{\partial^{2}u}{\partial x_{i}\partial
x_{j}}\frac{\partial^{2}\overline{v}}{\partial x_{i}\partial x_{j}}
dx+(\omega, z)_{0},\\
&&B((u,\omega),(v,z))=(\nabla(\frac{1}{n-1}u), \nabla v)_{0}+(\nabla
u, \nabla(\frac{n}{n-1}v))_{0}-(\frac{n}{n-1}\omega, v
)_{0}+(u,z)_{0},
\end{eqnarray*}
and if $(C2)$ holds, let
\begin{eqnarray*}
&&A((u,\omega),(v,z))=((\frac{1}{1-n}-\mu_{2})\Delta u, \Delta v)_{0}+(\mu_{2}\Delta u, \Delta v)_{0}+(\omega, z)_{0}\\
&&=((\frac{1}{1-n}-\mu_{2})\Delta u, \Delta
v)_{0}+\mu_{2}\int\limits_{\Omega} \sum\limits_{1\leq i,j\leq d}
\frac{\partial^{2}u}{\partial x_{i}\partial
x_{j}}\frac{\partial^{2}\overline{v}}{\partial x_{i}\partial x_{j}}
dx
+(\omega, z)_{0},\\
&&B((u,\omega),(v,z))=(\nabla(\frac{1}{1-n}u), \nabla v)_{0}+(\nabla
u, \nabla(\frac{n}{1-n}v))_{0}-(\frac{n}{1-n}\omega, v
)_{0}+(u,z)_{0},
\end{eqnarray*}
where $\mu_{1}>0$ and $\mu_{2}>0$ are chosen as  good approximations of $\min(\frac{1}{n-1})$ and $\min(\frac{1}{1-n})$ respectively such that $\frac{1}{n-1}-\mu_{1}\geq 0$ and
$\frac{1}{1-n}-\mu_{2}\geq 0$.\\
Let $\lambda=k^{2}$, then (\ref{s2.9})-(\ref{s2.10}) can be
rewritten as: Find $\lambda\in \mathbb{C}$, $(u,\omega)\in
\mathbf{H}\setminus \{0\}$ such that
\begin{eqnarray}\label{s2.11}
A((u,\omega),(v,z)) =\lambda B((u,\omega),(v,z)),~~~\forall (v,z)\in
\mathbf{H}.
\end{eqnarray}
 \indent Next we shall see
that the discrete bilinear form of (\ref{s2.11}) satisfies the
uniform $\mathbf{H}_{h}$-ellipticity automatically for many
non-conforming elements (see
(\ref{s2.22r})).\\

\indent Thus we get the following.\\
\indent{\bf Theorem 2.1.}~The weak formulations (\ref{s2.11}) and
(\ref{s2.6}) are equivalent.\\
\indent{\bf Proof.}~
 If $(k^{2}, u)$ is an eigenpair of (\ref{s2.6}), then together
with (\ref{s2.8}) we get that $(k^{2}, u, \omega)$ is an eigenpair
of (\ref{s2.9})-(\ref{s2.10}), thus it is an eigenpair of
(\ref{s2.11}). Conversely, if $(k^{2}, u, \omega)$ satisfies
(\ref{s2.11}), then $(k^{2}, u, \omega)$ also satisfies
(\ref{s2.9})-(\ref{s2.10}); from (\ref{s2.10}) we get
$\omega=k^{2}u$, and substituting it into (\ref{s2.9}) we get
(\ref{s2.6}). The above argument indicates that (\ref{s2.11}) and
(\ref{s2.6}) are equivalent.~~~$\Box$\\

 \indent For simplicity, in the next discussion we assume that $(C1)$
holds. And the argument is the same if $(C2)$ holds.\\
\indent  It is obvious that $A(\cdot,\cdot)$ is a selfadjoint,
continuous sesquilinear form on $\mathbf{H}\times \mathbf{H}$, and
\begin{eqnarray}\label{s2.12}
A((v,z),(v,z))\geq \mu_{1}|v|_{2}^{2}+\|z\|_{0}^{2}\gtrsim
\|(v,z)\|_{\mathbf{H}}^{2},
\end{eqnarray}
i.e., $A(\cdot,\cdot)$ is coercive.\\
\indent We use $A(\cdot,\cdot)$ and
$\|\cdot\|_{A}=A(\cdot,\cdot)^{\frac{1}{2}}$ as an inner product and
norm on $\mathbf{H}$, respectively.\\
 \indent Obviously, $k=0$ is not an
eigenvalue since
$A((u,\omega),(u,\omega))=0$ implies  $(u,\omega)=0$.\\
\indent When $n\in W^{1,\infty}(\Omega)$, a simple calculation shows
that
\begin{eqnarray}\label{s2.13}
&&|B((f,g), (v, z))|\nonumber\\
&&~~~=|(\nabla(\frac{1}{n-1}f), \nabla v)_{0}+(\nabla f,
\nabla(\frac{n}{n-1}v))_{0}- (\frac{n}{n-1}g, v
)_{0}+(f,z)_{0}|\nonumber\\
&&~~~\lesssim \|f\|_{1}\|v\|_{1}+ \|f\|_{1}\|v\|_{1}+
\|g\|_{-1}\|v\|_{1}+ \|f\|_{1}\|z\|_{-1}\nonumber\\
&&~~~\lesssim
(\|f\|_{1}+\|g\|_{-1})(\|v\|_{1}+\|z\|_{-1})\nonumber\\
&&~~~\lesssim
\|(f,g)\|_{\mathbf{H}_{1}}\|(v,z)\|_{\mathbf{H}_{1}},~~~\forall
(f,g), (v, z)\in \mathbf{H}_{1}.
\end{eqnarray}
\indent We can see from (\ref{s2.13}) that for any given $(f,g)\in
\mathbf{H}_{1}$, $B((f,g), (v, z))$ is a continuous linear form on
$\mathbf{H}$.\\
\indent The source problem associated with (\ref{s2.11}) is as
follows: Find $(\psi,\varphi)\in \mathbf{H}$ such that
\begin{eqnarray}\label{s2.14}
A((\psi,\varphi),(v,z)) =B((f,g),(v,z)),~~~\forall (v,z)\in
\mathbf{H}.
\end{eqnarray}
From Lax-Milgram theorem we know that (\ref{s2.14}) has one and only
one solution. Therefore, we define the corresponding solution
operators $T: \mathbf{H}_{1}\to \mathbf{H}$ by
\begin{eqnarray}\label{s2.15}
A(T(f,g),(v,z)) =B((f,g), (v,z)),~~~\forall (v,z)\in \mathbf{H}.
\end{eqnarray}
Then (\ref{s2.11}) has the equivalent operator form:
\begin{eqnarray}\label{s2.16}
T(u,\omega)=\frac{1}{\lambda}(u,\omega).
\end{eqnarray}

\indent{\bf Theorem 2.2.}~ Suppose $n\in W^{1,\infty}(\Omega)$, then
$T: \mathbf{H}\to \mathbf{H}$ is compact, and $T: \mathbf{H}_{1}\to
\mathbf{H}_{1}$ is compact.\\
\indent{\bf Proof.}~ Let $(v,z)=T(f,g)$ in (\ref{s2.14}), then from
(\ref{s2.12}) and (\ref{s2.13}) we have
\begin{eqnarray*}
\|T(f,g)\|_{\mathbf{H}}^{2}\lesssim A(T(f,g),T(f,g)) =B((f,g),
T(f,g))\lesssim
\|(f,g)\|_{\mathbf{H}_{1}}\|T(f,g)\|_{\mathbf{H_{1}}},
\end{eqnarray*}
thus
\begin{eqnarray}\label{s2.17}
\|T(f,g)\|_{\mathbf{H}} \lesssim \|(f,g)\|_{\mathbf{H}_{1}},
\end{eqnarray}
which implies that $T: \mathbf{H}_{1}\to \mathbf{H}$ is continuous.
Because of the compact embedding $\mathbf{H}\hookrightarrow
\mathbf{H}_{1}$, $T: \mathbf{H}\to \mathbf{H}$ is compact  and $T:
\mathbf{H}_{1}\to \mathbf{H}_{1}$ is compact.~~~$\Box$\\

\indent Consider the dual problem of (\ref{s2.11}): Find
$\lambda^{*}\in \mathbb{C}$, $(u^{*},\omega^{*})\in
\mathbf{H}\setminus \{0\}$ such that
\begin{eqnarray}\label{s2.18}
A((v,z), (u^{*},\omega^{*})) =\overline{\lambda^{*}} B((v,z),
(u^{*},\omega^{*})),~~~\forall (v,z)\in \mathbf{H}.
\end{eqnarray}
\indent The source problem associated with (\ref{s2.18}) is as
follows: Find $(\psi^{*},\varphi^{*})\in \mathbf{H}$ such that
\begin{eqnarray}\label{s2.19}
A((v,z), (\psi^{*},\varphi^{*})) =B((v,z), (f,g)),~~~\forall
(v,z)\in \mathbf{H}.
\end{eqnarray}
Define the corresponding solution operators $T^{*}:
\mathbf{H}_{1}\to \mathbf{H}$ by
\begin{eqnarray}\label{s2.20}
A((v,z), T^{*}(f,g)) =B((v,z), (f,g)),~~~\forall (v,z)\in
\mathbf{H}.
\end{eqnarray}
Then (\ref{s2.18}) has the equivalent operator form:
\begin{eqnarray}\label{s2.21}
T^{*}(u^{*},\omega^{*})=\lambda^{*-1}(u^{*},\omega^{*}).
\end{eqnarray}
\indent It can be proved that $T^{*}$ is the adjoint operator of $T$
in the sense of inner product $A(\cdot,\cdot)$. In fact, from
(\ref{s2.15}) and (\ref{s2.20}) we have
\begin{eqnarray*}
A(T(f,g),(v, z))=B((f,g),(v, z))=A((f,g),T^{*}(v, z)),~~~\forall
(f,g),(v, z)\in \mathbf{H}.
\end{eqnarray*}
\indent Note that since $T^{*}$ is the adjoint operator of $T$, the
primal and dual eigenvalues are connected via
$\lambda=\overline{\lambda^{*}}$.\\
\indent Let $\pi_{h}$ be a shape-regular mesh with size $h$. Let
$\mathbf{H}_{h}=S^{h}\times S^{h}\subset \mathbf{H}_{1}$ and
$\mathbf{H}_{h}\not\subset \mathbf{H}$
 be a
non-conforming finite element space; for example, $S^{h}\subset
H_{0}^{1}(\Omega)$ is the finite element space associated with one
of the Adini element, Morley-Zienkiewicz element, modified
Zienkiewicz
element, 12-parameter triangle plate element and 15-parameter triangle plate element et. al.\\
\indent Let
\begin{eqnarray*}
&&A_{h}((u_{h},\omega_{h}),(v,z))
=\sum\limits_{\kappa\in\pi_{h}}\int\limits_{\kappa}
\{(\frac{1}{n-1}-\mu_{1})\Delta u_{h} \Delta \overline{v}\\
&&~~~~~~+\mu_{1} \sum\limits_{1\leq i,j\leq d}
\frac{\partial^{2}u_{h}}{\partial x_{i}\partial
x_{j}}\frac{\partial^{2}\overline{v}}{\partial x_{i}\partial x_{j}}
\}dx+(\omega_{h}, z)_{0}.
\end{eqnarray*}
Denote
\begin{eqnarray*}
&&A_{h}((v,z),(v,z)) \equiv\|v\|_{h}^{2}+\|z\|_{0}^{2}\equiv
\|(v,z)\|_{h}^{2}.
\end{eqnarray*}
\indent For the finite element spaces mentioned above, from
\cite{ciarlet} and Lemma 5.4.3 of \cite{shi}, we know that
$A_{h}(\cdot,\cdot)$ satisfies  the uniform $\mathbf{H}_{h}$-ellipticity. \\
\begin{eqnarray}\label{s2.22r}
&&A_{h}((v,z),(v,z))\gtrsim
\sum\limits_{\kappa\in\pi_{h}}|v|_{2,\kappa}^{2}+\|z\|_{0}^{2},~~~\forall
(v,z)\in \mathbf{H}_{h}.
\end{eqnarray}
 Thus $\|(v,z)\|_{h}$ is a norm in $\mathbf{H}_{h}$, and the
generalized Poincare-Friedrichs inequality holds:
\begin{eqnarray}\label{s2.23}
\|(v,z)\|_{\mathbf{H}_{1}}\lesssim \|(v,z)\|_{h},~~~\forall (v,z)\in
\mathbf{H}_{h}.
\end{eqnarray}
\indent The non-conforming finite element approximation of
(\ref{s2.11}) is given by the following: Find $\lambda_{h}\in
\mathbb{C}$, $(u_{h},\omega_{h})\in \mathbf{H}_{h}\setminus \{0\}$
such that
\begin{eqnarray}\label{s2.24}
A_{h}((u_{h},\omega_{h}),(v,z)) =\lambda_{h}
B((u_{h},\omega_{h}),(v,z)),~~~\forall (v,z)\in \mathbf{H}_{h}.
\end{eqnarray}
\indent Consider the approximate source problem: Find
$(\psi_{h},\varphi_{h})\in \mathbf{H}_{h}$ such that
\begin{eqnarray}\label{s2.25}
A_{h}((\psi_{h},\varphi_{h}),(v,z)) =B((f,g), (v,z)),~~~\forall
(v,z)\in \mathbf{H}_{h}.
\end{eqnarray}
We introduce the corresponding solution operator: $T_{h}:
\mathbf{H}_{1}\to \mathbf{H}_{h}$:
\begin{eqnarray}\label{s2.26}
A_{h}(T_{h}(f,g),(v,z)) =B((f,g), (v,z)),~~~\forall (v,z)\in
\mathbf{H}_{h}.
\end{eqnarray}
Then (\ref{s2.24}) has the operator form:
\begin{eqnarray}\label{s2.27}
T_{h}(u_{h},\omega_{h})=\frac{1}{\lambda_{h}}(u_{h},\omega_{h}).
\end{eqnarray}
\indent The non-conforming finite element approximation of
(\ref{s2.18}) is given by:
 Find
$\lambda_{h}^{*} \in \mathbb{C}$, $(u_{h}^{*},\omega_{h}^{*})\in
\mathbf{H}_{h}\setminus \{0\}$ such that
\begin{eqnarray}\label{s2.28}
A_{h}((v,z), (u_{h}^{*},\omega_{h}^{*})) =\overline{\lambda_{h}^{*}}
B((v,z), (u_{h}^{*},\omega_{h}^{*})),~~~\forall (v,z)\in
\mathbf{H}_{h}.
\end{eqnarray}
Define the solution operator $T_{h}^{*}: \mathbf{H}_{1}\to
\mathbf{H}_{h}$ satisfying
\begin{eqnarray}\label{s2.29}
A_{h}((v,z), T_{h}^{*}(f,g))=B((v,z), (f,g)),~~~\forall~(v,z)\in
\mathbf{H}_{h}.
\end{eqnarray}
And (\ref{s2.28}) has the following equivalent operator form
\begin{eqnarray}\label{s2.30}
T_{h}^{*}(u_{h}^{*},\omega_{h}^{*})=\lambda_{h}^{*-1}(u_{h}^{*},\omega_{h}^{*}).
\end{eqnarray}
\indent It can be proved that $T_{h}^{*}$ is the adjoint operator of
$T_{h}$ in the sense of inner product $A_{h}(\cdot,\cdot)$. In fact,
from (\ref{s2.26}) and (\ref{s2.29}) we have
\begin{eqnarray*}
A_{h}(T_{h}(u,\omega),(v, z))=B((u,\omega),(v,
z))=A_{h}((u,\omega),T_{h}^{*}(v, z)),~~~\forall (u,\omega),(v,
z)\in \mathbf{H}_{h}.
\end{eqnarray*}
Hence, the primal and dual eigenvalues are connected via
$\lambda_{h}=\overline{\lambda_{h}^{*}}$.\\
\indent Denote
$$\mathbb{S}=(\frac{d}{2}, 2].$$
\indent Define interpolation operator $I_{h}^{1}:
H_{0}^{2}(\Omega)\cap W^{3,p}(\Omega)\to S^{h}$ ($p\in \mathbb{S}$),
and define $I_{h}^{2}: L^{2}(\Omega)\to S^{h}$ by
\begin{eqnarray*}
(\varphi-I_{h}^{2}\varphi, z)_{0}=0,~~~\forall z\in S^{h}.
\end{eqnarray*}
And let $I_{h}(\psi,\varphi)=(I_{h}^{1}\psi, I_{h}^{2}\varphi)$.\\
\indent For the finite element spaces mentioned above, when $\psi\in
W^{3,p}(\Omega)$ with $p\in \mathbb{S}$, the following estimates are
valid:
\begin{eqnarray}\label{s2.32}
&&\|I_{h}^{1}\psi-\psi\|_{h}\lesssim
h^{1+(\frac{1}{2}-\frac{1}{p})d}\|\psi\|_{3,p},\\\label{s2.33}
&&\|I_{h}^{1}\psi-\psi\|_{s}\lesssim
h^{3-s+(\frac{1}{2}-\frac{1}{p})d}\|\psi\|_{3,p},~~~s=0,1,
\end{eqnarray}
and when $\varphi\in H_{0}^{1}(\Omega)$
\begin{eqnarray}\label{s2.34}
&&\|I_{h}^{2}\varphi-\varphi\|_{0}=\inf\limits_{v\in
S^{h}}\|\varphi-v\|_{0}\lesssim h\|\varphi\|_{1},\\\label{s2.35}
&&\|I_{h}^{2}\varphi-\varphi\|_{-1}=\sup\limits_{v\in
H_{0}^{1}(\Omega)}\frac{(I_{h}^{2}\varphi-\varphi,
v-I^2_{h}v)_{0}}{\|v\|_{1}} \lesssim h^{2}\|\varphi\|_{1}.
\end{eqnarray}

\section{The consistency term and Strang lemma}

Let $(\psi,\varphi)$ and $(\psi^{*},\varphi^{*})$ be the solutions
of (\ref{s2.14}) and (\ref{s2.19}), respectively. Define the
consistency terms: For any $(v,z)\in \mathbf{H}_{h}+\mathbf{H}$,
\begin{eqnarray}\label{s3.1}
&&D_{h}((\psi,\varphi),(v,z))=B((f,g), (v,z))-A_{h}((\psi,\varphi),
(v,z)),\\\label{s3.2}
&&D_{h}^{*}((v,z),(\psi^{*},\varphi^{*}))=B((v,z),(f,g))-A_{h}((v,z),(\psi^{*},\varphi^{*})).
\end{eqnarray}
The following estimations of the consistency term play an crucial
role in our analysis.
\begin{eqnarray}\label{s3.3}
&&|D_{h}((\psi,\varphi),(v,z))| \lesssim
h^{1+(\frac{1}{2}-\frac{1}{p})d}\|\psi\|_{3,p}\|v\|_{h},\\\label{s3.4}
&&|D_{h}^{*}((v,z),(\psi^{*},\varphi^{*}))| \lesssim
h^{1+(\frac{1}{2}-\frac{1}{p})d}\|\psi^{*}\|_{3,p}\|v\|_{h}.
\end{eqnarray}
Next, we will prove the estimates (\ref{s3.3}) and (\ref{s3.4}) of the consistency term.\\
\indent It is well known that the following (C3) is valid for the
non-conforming finite elements mentioned in Section 2 except Adini
element (see Section 2.6 in
\cite{shi}).\\
\indent (C3)~If $F$ is the common face of element $\kappa$ and
$\kappa'$, then
\begin{eqnarray}\label{s3.5}
\int_{F}\nabla(v|_{\kappa})ds=\int_{F}\nabla(v|_{\kappa'})ds,~~~\forall
v\in S^{h};
\end{eqnarray}
if $F$ is a face of element $\kappa$ and $F\in
\partial\Omega$, then
\begin{eqnarray}\label{s3.6}
\int_{F}\nabla(v|_{\kappa})ds=0,~~~\forall v\in S^{h}.
\end{eqnarray}
\indent Define the face and element average interpolation operators
\begin{eqnarray*}
&&P_{F}^{0}f=\frac{1}{meas(F)}\int_{F}fds,~~~R_{F}^{0}f=f-P_{F}^{0}f,\\
&&P_{\kappa}^{0}f=\frac{1}{meas(\kappa)}\int_{\kappa}fdx,~~~R_{\kappa}^{0}f=f-P_{\kappa}^{0}f,
\end{eqnarray*}
where element $\kappa\in\pi_{h}$ and
  $F$ is   an arbitrary element face of $\pi_{h}$.\\
A simple calculation shows that for arbitrary constant $C_{0}$,
\begin{eqnarray}\label{s3.7}
\frac{1}{meas(F)}\int_{F}(f-P_{F}^{0}f)C_{0}ds=0.
\end{eqnarray}

\indent{\bf Theorem 3.1.}~ Suppose that $\psi, \psi^{*}\in
W^{3,p}(\Omega)$ $(p\in \mathbb{S})$, and (C3) is valid. Then for
any $(v,z) \in \mathbf{H}_{h}+\mathbf{H}$, (\ref{s3.3}) and
(\ref{s3.4}) hold.\\
\indent{\bf Proof.}~ For any $(v,z)\in C_{0}^{\infty}(\Omega)\times
L^{2}(\Omega)$, by the Green's formula we deduce
\begin{eqnarray}\label{s3.8}
B((f,g), (v,z))&=&A((\psi,\varphi),
(v,z))\nonumber\\
&=&\sum\limits_{\kappa\in\pi_{h}}\int\limits_{\kappa}
-\nabla(\frac{1}{n-1} \Delta \psi)\cdot
 \nabla \overline{v}+\varphi \overline{z} dx.
\end{eqnarray}
Since $C_{0}^{\infty}(\Omega)$ is dense in $H_{0}^{1}(\Omega)$, for
any $(v,z)\in \mathbf{H}_{h}+\mathbf{H}$ the above (\ref{s3.8})
holds. Thus
\begin{eqnarray}\label{s3.9}
&&D_{h}((\psi,\varphi),(v,z))=B((f,g), (v,z))-A_{h}((\psi,\varphi), (v,z))\nonumber\\
&&~~~=\sum\limits_{\kappa\in\pi_{h}}\int\limits_{\kappa}
-\nabla(\frac{1}{n-1} \Delta \psi)\cdot
 \nabla \overline{v}+\varphi \overline{z} dx
- \sum\limits_{\kappa\in\pi_{h}}\int\limits_{\kappa}
\{(\frac{1}{n-1}-\mu_{1})\Delta \psi \Delta \overline{v}\nonumber\\
&&~~~~~~+ \mu_{1}\sum\limits_{1\leq i,j\leq d}
\frac{\partial^{2}u}{\partial x_{i}\partial
x_{j}}\frac{\partial^{2}\overline{v}}{\partial x_{i}\partial x_{j}}
\} dx-(\varphi,z)_{0}\nonumber\\
&&~~~=-\sum\limits_{\kappa\in\pi_{h}}\int\limits_{\partial\kappa}
\frac{1}{n-1} \Delta \psi
 \nabla \overline{v}\cdot\gamma ds\nonumber\\
&&~~~~~~+\mu_{1}\sum\limits_{\kappa\in\pi_{h}}\int\limits_{\kappa}
\{\Delta \psi \Delta \overline{v}- \sum\limits_{1\leq i,j\leq d}
\frac{\partial^{2}u}{\partial x_{i}\partial
x_{j}}\frac{\partial^{2}\overline{v}}{\partial x_{i}\partial x_{j}}
\}dx\nonumber\\
&&~~~=-\sum\limits_{\kappa\in\pi_{h}}\int\limits_{\partial\kappa}
\frac{1}{n-1} \Delta \psi
 \nabla \overline{v}\cdot\gamma ds\nonumber\\
&&~~~~~~+\mu_{1}\sum\limits_{\kappa\in\pi_{h}}\int\limits_{\kappa}
(\sum\limits_{1\leq i\not=j\leq d} \frac{\partial^{2}\psi}{\partial
x_{i}^{2}}\frac{\partial^{2} \overline{v}}{\partial x_{j}^{2}} -
\sum\limits_{1\leq i\not=j\leq d}\frac{\partial^{2}\psi}{\partial
x_{i}\partial x_{j}}\frac{\partial^{2}\overline{v}}{\partial
x_{i}\partial
x_{j}})dx\nonumber\\
&&~~~=-\sum\limits_{\kappa\in\pi_{h}}\int\limits_{\partial\kappa}
\frac{1}{n-1} \Delta \psi
 \nabla \overline{v}\cdot\gamma ds
 +\mu_{1}\sum\limits_{\kappa\in\pi_{h}}\int\limits_{\partial\kappa}
\sum\limits_{1\leq i\not=j\leq d} \frac{\partial^{2}\psi}{\partial
x_{i}^{2}}\frac{\partial \overline{v}}{\partial x_{j}}\gamma_{j}
ds\nonumber\\
&&~~~~~~-\mu_{1}\sum\limits_{\kappa\in\pi_{h}}\int\limits_{\partial\kappa}
\sum\limits_{1\leq i\not=j\leq d}\frac{\partial^{2}\psi}{\partial
x_{i}\partial x_{j}}\frac{\partial\overline{v}}{\partial
x_{j}}\gamma_{i}
ds\nonumber\\
&&~~~\equiv I_{1}+I_{2}+I_{3}.
\end{eqnarray}
\indent Note that (C3) and (\ref{s3.7}) are valid, and for all $
v\in H_{0}^{2}(\Omega)$ (\ref{s3.5})-(\ref{s3.6}) also  hold, we
deduce that $\forall (v,z) \in \mathbf{H}_{h}+\mathbf{H}$,
\begin{eqnarray}\label{s3.10}
I_{1}=-\sum\limits_{\kappa\in\pi_{h}}\sum\limits_{F\in\partial\kappa}\int\limits_{F}
R_{F}^{0}(\frac{1}{n-1} \Delta \psi) R_{F}^{0}( \nabla
\overline{v}\cdot\gamma )ds.
\end{eqnarray}
Let $\hat{\kappa}$ is a reference element, $\kappa$ and
$\hat{\kappa}$ be affine-equivalent. When $\hat{w}\in
W_{1,\iota}(\hat{\kappa})$ and $1\leq
\rho<\frac{(d-1)\iota}{d-\iota}$, by the trace theorem we get
$W_{1,\iota}(\hat{\kappa})\hookrightarrow
L_{\rho}(\partial\hat{\kappa})$, thus we deduce the following trace
inequality:
\begin{eqnarray}\label{s3.11}
&&\int\limits_{\partial\kappa}
|w|^{\rho}ds=\int\limits_{\partial\hat{\kappa}}
|\hat{w}|^{\rho}\frac{|\partial\kappa|}{\partial\hat{\kappa}}d\hat{s}\lesssim
h_{\kappa}^{d-1}\|\hat{w}\|_{0,\rho,\partial\hat{\kappa}}^{\rho}
\lesssim h_{\kappa}^{d-1}\|\hat{w}\|_{1,\iota,\hat{\kappa}}^{\rho}\nonumber\\
&&~~~\lesssim
h_{\kappa}^{d-1}(\|\hat{w}\|_{0,\iota,\hat{\kappa}}^{\rho}+|\hat{w}|_{1,\iota,\hat{\kappa}}^{\rho})
\lesssim h_{\kappa}^{d-1}(h_{\kappa}^{-\frac{\rho
d}{\iota}}\|w\|_{0,\iota,\kappa}^{\rho}
+h_{\kappa}^{\rho-\frac{\rho d}{\iota}}|w|_{1,\iota,\kappa}^{\rho})\nonumber\\
&&~~~\lesssim h_{\kappa}^{d-\frac{\rho
d}{\iota}-1}\|w\|_{0,\iota,\kappa}^{\rho}+h_{\kappa}^{\rho+d-\frac{\rho
d}{\iota}-1}|w|_{1,\iota,\kappa}^{\rho},~~~\forall \kappa\in
\pi_{h}.
\end{eqnarray}
Since $p\in \mathbb{S}$, $W_{1,p}(\hat{\kappa})\hookrightarrow
L_{\rho}(\partial\hat{\kappa})$ with $\rho\in(d-1,
\frac{(d-1)p}{d-p})$. Choose $\frac{1}{\rho'}=1-\frac{1}{\rho}$,
then $\rho'<\frac{d-1}{d-2}$ and
$W_{1,2}(\hat{\kappa})\hookrightarrow
L_{\rho'}(\partial\hat{\kappa})$. And thus, by the H$\ddot{o}$lder
inequality, the trace inequality (\ref{s3.11}) and the interpolation
error estimate we deduce that
\begin{eqnarray}\label{s3.12}
&& |I_{1}|\lesssim
\sum\limits_{\kappa\in\pi_{h}}\sum\limits_{F\in\partial\kappa}
\|R_{F}^{0}(\frac{1}{n-1} \Delta \psi)\|_{0,\rho,F}
\|R_{F}^{0}( \nabla \overline{v}\cdot\gamma )\|_{0,\rho',F}\nonumber\\
&&~~~\lesssim
\sum\limits_{\kappa\in\pi_{h}}\sum\limits_{F\in\partial\kappa}
  \|R_{\kappa}^{0}(\frac{1}{n-1} \Delta
\psi)\|_{0,\rho,F} \|R_{\kappa}^{0}( \nabla \overline{v}\cdot\gamma
)\|_{0,\rho',F}\nonumber\\
&&~~~\lesssim \sum\limits_{\kappa\in\pi_{h}}
(h_{\kappa}^{d-\frac{\rho d}{p}-1}\|R_{\kappa}^{0}(\frac{1}{n-1}
\Delta \psi)\|_{0,p,\kappa}^{\rho} +h_{\kappa}^{\rho+d-\frac{\rho
d}{p}-1}|R_{\kappa}^{0}(\frac{1}{n-1} \Delta
\psi)|_{1,p,\kappa}^{\rho})^{\frac{1}{\rho}}\nonumber\\
&&~~~~~~\times(h_{\kappa}^{d-\frac{\rho'd}{2}-1}\|R_{\kappa}^{0}(
\nabla v\cdot\gamma )\|_{0,\kappa}^{\rho'}
+h_{\kappa}^{\rho'+d-\frac{\rho'd}{2}-1}|R_{\kappa}^{0}( \nabla
\overline{v}\cdot\gamma
)|_{1,\kappa}^{\rho'})^{\frac{1}{\rho'}}\nonumber\\
&&~~~\lesssim \sum\limits_{\kappa\in\pi_{h}}
(h_{\kappa}^{\rho+d-\frac{\rho
d}{p}-1})^{\frac{1}{\rho}}\|\psi\|_{3,p,\kappa}
\times(h_{\kappa}^{\rho'+d-\frac{\rho'd}{2}-1})^{\frac{1}{\rho'}}\|v\|_{2,\kappa}\nonumber\\
&&~~~\lesssim
h^{1+(\frac{1}{2}-\frac{1}{p})d}\|\psi\|_{3,p}\|v\|_{h},~~~\forall
(v,z) \in \mathbf{H}_{h}+\mathbf{H}.
\end{eqnarray}
Similarly we deduce
\begin{eqnarray}\label{s3.13}
&& |I_{2}|\lesssim
h^{1+(\frac{1}{2}-\frac{1}{p})d}\|\psi\|_{3,p}\|v\|_{h},~~~\forall
(v,z) \in \mathbf{H}_{h}+\mathbf{H},\\\label{s3.14} &&
|I_{3}|\lesssim
h^{1+(\frac{1}{2}-\frac{1}{p})d}\|\psi\|_{3,p}\|v\|_{h},~~~\forall
(v,z) \in \mathbf{H}_{h}+\mathbf{H}.
\end{eqnarray}
Substituting (\ref{s3.12}), (\ref{s3.13}) and (\ref{s3.14}) into
(\ref{s3.9}) we get (\ref{s3.3}).\\
\indent Using the same argument as above,  we can prove
(\ref{s3.4}).~~~$\Box$\\

\indent Next, we shall analyze Adini rectangle element
approximation. We suppose that $\Omega\subset \mathbb{R}^{2}$, and
the boundary of $\Omega$ and the edges of elements are parallel to
the coordinate axis.
   Although
(C3) is not valid,  thanks to \cite{ciarlet}, we can prove
(\ref{s3.3}) and (\ref{s3.4}) still hold.\\

\indent{\bf Theorem 3.2.}~ Suppose that $S^{h}$ is Adini element
space, then for any $ (v,z) \in \mathbf{H}_{h}$ (\ref{s3.3}) and
(\ref{s3.4}) are valid.\\
\indent{\bf Proof.}~ We shall analyze the terms $I_{1}$, $I_{2}$ and
$I_{3}$ on the right-hand side of (\ref{s3.9}). Noticing that the
edges of elements are parallel to the coordinate axis, using the
proof method of Theorem 50.1 in \cite{ciarlet}, we can deduce that
for any $(v,z)\in \mathbf{H}_{h}$,
\begin{eqnarray*}
|I_{1}|+|I_{2}|\lesssim
h^{1+(\frac{1}{2}-\frac{1}{p})d}\|\psi\|_{3,p}\|v\|_{h}.
\end{eqnarray*}
And from line 11 on page 304 in \cite{ciarlet}, we see
\begin{eqnarray*}
I_{3}=0.
\end{eqnarray*}
Substituting the above estimates into (\ref{s3.9}) we get
(\ref{s3.3}). Similarly we can prove (\ref{s3.4}).~~~$\Box$\\

\indent The following lemma is a generalization of Strang Lemma
(1972).\\

\indent{\bf Lemma 3.1.}~
 Let $(\psi,\varphi)$ be the solution of
(\ref{s2.14}) and $(\psi_{h},\varphi_{h})$ be the solution of
(\ref{s2.25}), then
\begin{eqnarray}\label{s3.15}
&&\inf\limits_{(v,z)\in
\mathbf{H}_{h}}\|(\psi,\varphi)-(v,z)\|_{h}+\sup\limits_{(v,z)\in
\mathbf{H}_{h}\setminus\{0\}}\frac{
D_{h}((\psi,\varphi),(v,z))}{\|(v,z)\|_{h}}\nonumber\\
&&~~~\lesssim
\|(\psi,\varphi)-(\psi_{h},\varphi_{h})\|_{h}\nonumber\\
&&~~~\lesssim \inf\limits_{(v,z)\in
\mathbf{H}_{h}}\|(\psi,\varphi)-(v,z)\|_{h}+\sup\limits_{(v,z)\in
\mathbf{H}_{h}\setminus\{0\}} \frac{
D_{h}((\psi,\varphi),(v,z))}{\|(v,z)\|_{h}}.
\end{eqnarray}
Let $(\psi^{*},\varphi^{*})$ be the solution of (\ref{s2.19}) and
$(\psi_{h}^{*},\varphi_{h}^{*})$ be its finite element solution,
then
\begin{eqnarray}\label{s3.16}
&&\inf\limits_{(v,z)\in
\mathbf{H}_{h}}\|(\psi^{*},\varphi^{*})-(v,z)\|_{h}+\sup\limits_{(v,z)\in
\mathbf{H}_{h}\setminus\{0\}}\frac{
D_{h}((v,z), (\psi^{*},\varphi^{*}))}{\|(v,z)\|_{h}}\nonumber\\
&&~~~\lesssim
\|(\psi^{*},\varphi^{*})-(\psi_{h}^{*},\varphi_{h}^{*})\|_{h}\nonumber\\
&&~~~\lesssim \inf\limits_{(v,z)\in
\mathbf{H}_{h}}\|(\psi^{*},\varphi^{*})-(v,z)\|_{h}+\sup\limits_{(v,z)\in
\mathbf{H}_{h}\setminus\{0\}} \frac{ D_{h}((v,z),
(\psi^{*},\varphi^{*}))}{\|(v,z)\|_{h}}.
\end{eqnarray}
\indent{\bf Proof.}~ For any $(v,z)\in \mathbf{H}_{h}$,
\begin{eqnarray*}
&&\|(\psi_{h},\varphi_{h})-(v,z)\|_{h}^{2}= A_{h}((\psi_{h},\varphi_{h})-(v,z),(\psi_{h},\varphi_{h})-(v,z))\\
&&~~~=A_{h}((\psi,\varphi)-(v,z),(\psi_{h},\varphi_{h})-(v,z))+B((f,g),(\psi_{h},\varphi_{h})-(v,z))\nonumber\\
&&~~~~~~-A_{h}((\psi,\varphi),(\psi_{h},\varphi_{h})-(v,z)).
\end{eqnarray*}
When $\|(\psi_{h},\varphi_{h})-(v,z)\|_{h}\not= 0$, dividing it in
both sides of the above we obtain
\begin{eqnarray*}
&&\|(\psi_{h},\varphi_{h})-(v,z)\|_{h} \leq
\|(\psi,\varphi)-(v,z)\|_{h}\nonumber\\
&&~~~~~~-\frac{ A_{h}((\psi,\varphi),(\psi_{h},\varphi_{h})-(v,z))
-B((f,g),(\psi_{h},\varphi_{h})-(v,z))}{\|(\psi_{h},\varphi_{h})-(v,z)\|_{h}}\\
&&~~~\lesssim \|(\psi,\varphi)-(v,z)\|_{h}+\sup\limits_{(v,z)\in
\mathbf{H}_{h}\setminus\{0\}}
\frac{D_{h}((\psi,\varphi),(v,z))}{\|(v,z)\|_{h}}.
\end{eqnarray*}
This together with the triangular inequality
\begin{eqnarray*}
\|(\psi,\varphi)-(\psi_{h},\varphi_{h})\|_{h}\leq
\|(\psi,\varphi)-(v,z)\|_{h}+\|(v,z)-(\psi_{h},\varphi_{h})\|_{h}
\end{eqnarray*}
yields the second inequality of (\ref{s3.15}). From
\begin{eqnarray*}
A_{h}((\psi,\varphi)-(\psi_{h},\varphi_{h}),(v,z))\leq
\|(\psi,\varphi)-(\psi_{h},\varphi_{h})\|_{h}\|(v,z)\|_{h},~~~~\forall~
(v,z)\in S^{h},
\end{eqnarray*}
we get
\begin{eqnarray*}
\|(\psi,\varphi)-(\psi_{h},\varphi_{h})\|_{h}\geq \frac{
A_{h}((\psi,\varphi),(v,z))-A_{h}((\psi_{h},\varphi_{h}),(v,z))}{\|(v,z)\|_{h}}=-\frac{
D_{h}((\psi,\varphi),(v,z))}{\|(v,z)\|_{h}},
\end{eqnarray*}
which together with
$\|(\psi,\varphi)-(\psi_{h},\varphi_{h})\|_{h}\geq
\inf\limits_{(v,z)\in S^{h}}\|(\psi,\varphi)-(v,z)\|_{h}$ we obtain
the first inequality of (\ref{s3.15}).\\
\indent Similarly we can prove (\ref{s3.16}). The proof is
completed.~~~$\Box$\\

\indent By Lemma 3.1, we get:\\
\indent{\bf Theorem 3.3.}~ Suppose that $\psi, \psi^{*}\in
W^{3,p}(\Omega)$ $(p\in \mathbb{S})$, for any $(v,z) \in
\mathbf{H}_{h}$ (\ref{s3.3}) and (\ref{s3.4}) hold. Then
\begin{eqnarray}\label{s3.17}
&&\|(\psi,\varphi)-(\psi_{h},\varphi_{h})\|_{h}\lesssim
h^{1+(\frac{1}{2}-\frac{1}{p})d}(\|\psi\|_{3,p}+\|\varphi\|_{1}),\\\label{s3.18}
&&\|(\psi^{*},\varphi^{*})-(\psi_{h}^{*},\varphi_{h}^{*})\|_{h}
\lesssim
h^{1+(\frac{1}{2}-\frac{1}{p})d}(\|\psi^{*}\|_{3,p}+\|\varphi^{*}\|_{1}).
\end{eqnarray}
\indent{\bf Proof.}~
 By the interpolation error estimates (\ref{s2.32}) and
(\ref{s2.34}), we get
\begin{eqnarray}\label{s3.19}
&&\inf\limits_{(v,z)\in
\mathbf{H}_{h}}\|(\psi,\varphi)-(v,z)\|_{h}\lesssim
\|(\psi,\varphi)-I_{h}(\psi,\varphi)\|_{h}\nonumber\\
&&~~~=\|\psi-I_{h}^{1}\psi\|_{h}
+\|\varphi-I_{h}^{2}\varphi\|_{0}\nonumber\\
&&~~~\lesssim
h^{1+(\frac{1}{2}-\frac{1}{p})d}\|\psi\|_{3,p}+h\|\varphi\|_{1}\lesssim
h^{1+(\frac{1}{2}-\frac{1}{p})d}(\|\psi\|_{3,p}+\|\varphi\|_{1}).
\end{eqnarray}
Substituting (\ref{s3.19}) and (\ref{s3.3}) into (\ref{s3.15}) we
get (\ref{s3.17}). By the same argument we can prove (\ref{s3.18}).
The proof is completed.~~~$\Box$\\

\indent{\bf Remark 3.1.}~~{\em We tried to use the Nitsche technique
to prove the error estimate in $\|\cdot\|_{\mathbf{H}_{1}}$ is of
higher order than that in $\|\cdot\|_{h}$, but failed because of the
non-symmetry of right-hand sides that involves derivatives, of (\ref{s2.14})
and (\ref{s2.19}).}

\section{The error analysis of the non-conforming element eigenvalues}
\indent Let $(\lambda,u,\omega)$ and
$(\lambda^{*},u^{*},\omega^{*})$
 be the eigenpair of (\ref{s2.11}) and (\ref{s2.18}), respectively. Then from (\ref{s3.1}) and
 (\ref{s3.2}) we get that for any $(v,z)\in
\mathbf{H}_{h}+\mathbf{H}$,
 \begin{eqnarray}\label{s4.1}
&&D_{h}((u,\omega),(v,z))=B(\lambda(u,\omega),
(v,z))-A_{h}((u,\omega), (v,z)),\\\label{s4.2}
&&D_{h}^{*}((v,z),(u^{*},\omega^{*}))=B((v,z),\lambda^{*}(u^{*},\omega^{*}))-A_{h}((v,z),(u^{*},\omega^{*})).
\end{eqnarray}
\indent We need the following regularity assumption:\\
\indent {\bf $R(\Omega)$.}~~{\em For any $\xi\in H^{-1}(\Omega)$,
there exists $\psi\in W_{3,p_{0}}(\Omega)$ satisfying
\begin{eqnarray*}
\Delta(\frac{1}{n-1}\Delta
\psi)=\xi,~~~in~\Omega,~~~\psi=\frac{\partial \psi}{\partial
\nu}=0~~~ on~\partial \Omega,
\end{eqnarray*}
and
\begin{eqnarray}\label{s4.3}
\|\psi\|_{3,p_{0}}\leq C_{R} \|\xi\|_{-1},
\end{eqnarray}
where $p_{0}\in \mathbb{S}$, $C_{R}$ denotes the  prior constant
dependent on the $n(x)$ and $\Omega$ but independent of the
right-hand
side $\xi$ of the equation.}\\
 \indent It is well known that (\ref{s4.3}) is valid when $n$ and $\partial \Omega$ are appropriately smooth.
For example, when $\Omega\subset \mathbb{R}^{2}$ is a convex
polygon, from Theorem 2 in \cite{blum}, we can get that
$p_{0}=2$.\\
\indent Consider the source problem associated with (\ref{s2.5}) and
(\ref{s2.7}):
\begin{eqnarray}\label{s4.4}
&&\Delta(\frac{1}{n-1}\Delta \psi)=-\frac{n}{n-1}\Delta f-\Delta
(\frac{1}{n-1}f)-\frac{n}{n-1}g,\\\label{s4.5} &&\varphi=f.
\end{eqnarray}
When $n\in W^{1,\infty}(\Omega)\cap H^{2}(\Omega)$ and $f$ is
appropriately smooth, from $R(\Omega)$ we can deduce that $\psi\in
W^{3,p_{0}}(\Omega)$ and
\begin{eqnarray}\label{s4.6}
&&\|\psi\|_{3,p_{0}}\leq C_{R}\|-\frac{n}{n-1}\Delta f-\Delta
(\frac{1}{n-1}f)-\frac{n}{n-1}g\|_{-1}\nonumber\\
&&~~~\lesssim
\|f\|_{1}+\|g\|_{-1}\lesssim\|(f,g)\|_{\mathbf{H}_{1}}\\\label{s4.7}
&&\|\varphi\|_{1}=\|f\|_{1}.
\end{eqnarray}
\indent In this paper, for simplicity, we assume that the dual and
primal problems have the same regularity.\\

\indent{\bf Theorem 4.1.}~ Assume $n\in W^{1,\infty}(\Omega)\cap
H^{2}(\Omega)$, (\ref{s2.2}) and $R(\Omega)$ hold, and for any
$(v,z) \in \mathbf{H}_{h}$ (\ref{s3.3})-(\ref{s3.4}) are valid. Then
\begin{eqnarray}\label{s4.8}
&&\|T-T_{h}\|_{\mathbf{H}_{1}}\lesssim
h^{1+(\frac{1}{2}-\frac{1}{p_{0}})d},\\\label{s4.9}
&&\|T^{*}-T_{h}^{*}\|_{\mathbf{H}_{1}}\lesssim
h^{1+(\frac{1}{2}-\frac{1}{p_{0}})d}.
\end{eqnarray}
\indent{\bf Proof.}~
 For any $(f,g)\in \mathbf{H}_{1}$, with
$\|(f,g)\|_{\mathbf{H}_{1}}=1$, there is $(f_{j},g)\in
C_{0}^{\infty}(\Omega)\times H^{-1}(\Omega)$, such that
$$\|(f,g)-(f_{j},g)\|_{\mathbf{H}_{1}}=\|f-f_{j}\|_{1}\leq h.$$
By (\ref{s2.23}), (\ref{s2.32})-(\ref{s2.35}), (\ref{s3.17}) and
(\ref{s4.6}) we deduce
\begin{eqnarray*}
&&\|(T-T_{h})(f_{j},g)\|_{\mathbf{H}_{1}}\leq
\|T(f_{j},g)-I_{h}T(f_{j},g)\|_{\mathbf{H}_{1}}+\|I_{h}T(f_{j},g)-T_{h}(f_{j},g)\|_{\mathbf{H}_{1}}\\
&&~~~\leq
\|T(f_{j},g)-I_{h}T(f_{j},g)\|_{\mathbf{H}_{1}}+\|I_{h}T(f_{j},g)-T_{h}(f_{j},g)\|_{h}\\
&&~~~\leq
Ch^{1+(\frac{1}{2}-\frac{1}{p_{0}})d}\|T(f_{j},g)\|_{W^{3,p_{0}}(\Omega)\times
H^{1}(\Omega)}
+\|T(f_{j},g)-T_{h}(f_{j},g)\|_{h}\\
&&~~~\leq
h^{1+(\frac{1}{2}-\frac{1}{p_{0}})d}\|T(f_{j},g)\|_{W^{3,p_{0}}(\Omega)\times
H^{1}(\Omega)}\leq
Ch^{1+(\frac{1}{2}-\frac{1}{p_{0}})d}\|(f_{j},g)\|_{\mathbf{H}_{1}}.
\end{eqnarray*}
From (\ref{s2.26}) we know that $T_{h}$ has a upper bound uniformly
with respect to $h$. Thus we have
\begin{eqnarray*}
&&\|(T-T_{h})(f,g)\|_{\mathbf{H}_{1}}\leq
\|(T-T_{h})((f,g)-(f_{j},g))\|_{\mathbf{H}_{1}}+\|(T-T_{h})(f_{j},g)\|_{\mathbf{H}_{1}}\\
&&~~~\leq
(\|T\|_{\mathbf{H}_{1}}+\|T_{h}\|_{\mathbf{H}_{1}})\|(f,g)-(f_{j},g)\|_{\mathbf{H}_{1}}+Ch^{1+(\frac{1}{2}-\frac{1}{p_{0}})d}\|(f_{j},g)\|_{\mathbf{H}_{1}}\\
&&~~~\leq
(\|T\|_{\mathbf{H}_{1}}+\|T_{h}\|_{\mathbf{H}_{1}})h+Ch^{1+(\frac{1}{2}-\frac{1}{p_{0}})d}(\|(f_{j},g)-(f,g)\|_{\mathbf{H}_{1}}+\|(f,g)\|_{\mathbf{H}_{1}})\\
&&~~~\leq
(\|T\|_{\mathbf{H}_{1}}+\|T_{h}\|_{\mathbf{H}_{1}}+C)h^{1+(\frac{1}{2}-\frac{1}{p_{0}})d}.
\end{eqnarray*}
And by the definition of operator norm we have
\begin{eqnarray*}
\|T-T_{h}\|_{\mathbf{H}_{1}} =\sup\limits_{(f,g)\in \mathbf{H}_{1},
\|(f,g)\|_{\mathbf{H}_{1}}=1}\|(T-T_{h})(f,g)\|_{\mathbf{H}_{1}}\lesssim
h^{1+(\frac{1}{2}-\frac{1}{p_{0}})d}.
\end{eqnarray*}
Hence, (\ref{s4.8}) is valid. Similarly we can deduce (\ref{s4.9}).
The proof is completed.~~~$\Box$\\

\indent In this paper, we suppose that $\lambda$ be an eigenvalue of
(\ref{s2.11}) with the algebraic multiplicity $q$ and the ascent
$\alpha$. Then $\lambda^*=\overline{\lambda}$ is an eigenvalue of
(\ref{s2.18}). Since $\|T_{h}-T\|_{\mathbf{H}_{1}}\to 0$, $q$
eigenvalues $\lambda_{1,h},\cdots,\lambda_{q,h}$ of (\ref{s2.24})
will converge to $\lambda$.\\
\indent Let $E$ be the spectral projection associated with $T$ and
$\lambda$, then $R(E)=N((\lambda^{-1}-T))$ is the space of
generalized eigenfunctions associated with $\lambda$ and $T$, where
$R$ denotes the range and $N$ denotes the null space. Let $E_{h}$ be
the spectral projection associated with $T_{h}$ and the eigenvalues
$\lambda_{1,h},\cdots,\lambda_{q,h}$, then $ R(E_{h}) $ is the space
spanned by all generalized eigenfunctions corresponding to all
eigenvalues $\lambda_{1,h},\cdots,\lambda_{q,h}$. In view of the
adjoint problem (\ref{s2.18}) and (\ref{s2.28}), the definitions of
$E^{*}$, $R(E^{*})$, $E_{h}^{*}$ and $R(E_{h}^{*})$ are analogous
to $E$, $R(E)$, $E_{h}$ and $R(E_{h})$  (see \cite{babuska}).\\
\indent Let $\lambda_{h}\in\{\lambda_{1,h},\cdots,\lambda_{q,h}\}$.
From \cite{babuska} we get the following results.

\indent{\bf Theorem 4.2.}~ Assume that the conditions of Theorem 4.1
are valid. Let $(u_{h},\omega_{h})$ be eigenfunction corresponding
to $\lambda_{h}$ and $\|(u_{h},\omega_{h})\|_h=1$. Then there exists
eigenfunction $(u,\omega)$ corresponding to $\lambda$, such that
\begin{eqnarray}\label{s4.10}
&&\|(u_{h},\omega_{h})-(u,\omega)\|_{\mathbf{H}_{1}}\lesssim
\|(T-T_{h})|_{R(E)}\|_{\mathbf{H}_{1}}^{\frac{1}{\alpha}},\\
\label{s4.11} &&|\lambda_{h}-\lambda|\lesssim
\|(T-T_{h})|_{R(E)}\|_{\mathbf{H}_{1}}^{\frac{1}{\alpha}}.\\\label{s4.12}
&&|(\frac{1}{q}\sum\limits_{i=1}^{q}\lambda_{i,h}^{-1})^{-1}-\lambda|\lesssim
\|(T-T_{h})|_{R(E)}\|_{\mathbf{H}_{1}}.
\end{eqnarray}
Furthermore assume $R(E)\subset W^{3,p}(\Omega)$ $(p\in
\mathbb{S})$, then
\begin{eqnarray}\label{s4.13}
\|(T-T_{h})|_{R(E)}\|_{\mathbf{H}_{1}}\lesssim
h^{1+(\frac{1}{2}-\frac{1}{p})d},
\end{eqnarray}
and
\begin{eqnarray}\label{s4.14}
&&\|(u_{h},\omega_{h})-(u,\omega)\|_{\mathbf{H}_{1}} \lesssim
h^{\frac{1}{\alpha}+(\frac{1}{2}-\frac{1}{p})\frac{d}{\alpha}},\\\label{s4.15}
&&\|(u_{h},\omega_{h})-(u,\omega)\|_{h}\lesssim
h^{\frac{1}{\alpha}+(\frac{1}{2}-\frac{1}{p})\frac{d}{\alpha}},
\end{eqnarray}
with $\|(u,\omega)\|_h=1$.\\
\indent{\bf Proof.}~ From Theorem 4.1 we know
$\|T-T_{h}\|_{\mathbf{H}_{1}}\to 0~(h\to
 0)$, thus from Theorem 7.4, Theorem 7.3 and Theorem
 7.2
of \cite{babuska} we get (\ref{s4.10}), (\ref{s4.11}) and
(\ref{s4.12}), respectively. By the way to show (\ref{s4.8}), we get
(\ref{s4.13}). Substituting (\ref{s4.13}) into (\ref{s4.10}), we get
(\ref{s4.14}). By calculation we get
\begin{eqnarray}\label{s4.16}
&&\|(u_{h},\omega_{h})-(u,\omega)\|_{h}
=\|\lambda_{h}T_{h}(u_{h},\omega_{h})-\lambda T(u,\omega)\|_{h}\nonumber\\
&&~~~\leq \|\lambda_{h}T_{h}(u_{h},\omega_{h})-\lambda
T_{h}(u,\omega)\|_{h} +\|\lambda T_{h}(u,\omega)
-\lambda T(u,\omega)\|_{h}\nonumber\\
&&~~~\lesssim \|\lambda T_{h}(u,\omega) -\lambda T(u,\omega)\|_{h}+
\|\lambda_{h}(u_{h},\omega_{h})-\lambda
(u,\omega)\|_{\mathbf{H}_{1}}\nonumber\\
&&~~~\lesssim |\lambda|\|T(u,\omega)-T_{h}(u,\omega)\|_{h}+
h^{\frac{1}{\alpha}+(\frac{1}{2}-\frac{1}{p})\frac{d}{\alpha}}.
\end{eqnarray}
Combining (\ref{s3.17}) with the above relation we get
(\ref{s4.15}). By calculation we have
\begin{eqnarray}\label{s4.17}
&&\|(u_{h},\omega_{h})-\frac{(u,\omega)}{\|(u,\omega)\|_{h}}\|_{\mathbf{s}}\lesssim
\|(u_{h},\omega_{h})-(u,\omega)\|_{h}\nonumber\\
&&~~~~~~+\|(u_{h},\omega_{h})-(u,\omega)\|_{\mathbf{s}},
~~~\mathbf{s}=\mathbf{H}_{1}, h,
\end{eqnarray}
thus, when replacing $(u,\omega)$ by
$\frac{(u,\omega)}{\|(u,\omega)\|_{h}}$, (\ref{s4.14}) and
(\ref{s4.15}) also hold.~~~$\Box$\\

\indent
 Starting from (\ref{s4.11}), if we use $\|(T-T_{h})|_{R(E)}\|_{\mathbf{H}_{1}}$
 we can not derive the optimal estimates for the eigenvalue when the eigenfunction is smooth on concave domain
 because the error estimate in $\|\cdot\|_{\mathbf{H}_{1}}$ depends on
 the Nitsche technique and
the regularity.
 To avoid this problem, we employ a new method and give an identity in the following lemma, and use it
to prove the optimal error estimates of non-conforming element
eigenvalues. The identity and proof method are also valid for
general nonselfadjoint eigenvalue problems.\\

\indent{\bf Lemma 4.1.}~ Let $(\lambda,u,\omega)$ and
$(\lambda^{*},u^{*},\omega^{*})$ be the eigenpairs of (\ref{s2.11})
and (\ref{s2.18}) respectively.
 Then for any $(v,z), (v^{*},z^{*})\in
 \mathbf{H}_{h}$, when
 $B((v,z), (v^{*},z^{*}))\not=0$ it is valid that
\begin{eqnarray}\label{s4.18}
&&\frac{A_{h}((v,z), (v^{*},z^{*}))}{B((v,z),
(v^{*},z^{*}))}-\lambda =\frac{A_{h}((u,\omega)-(v,z),
(u^{*},\omega^{*})-(v^{*},z^{*}))}{B((v,z),(v^{*},z^{*}))}\nonumber\\
&&~~~~~~-\lambda \frac{B((u,\omega)-(v,z),
(u^{*},\omega^{*})-(v^{*},z^{*}))}{B((v,z),(v^{*},z^{*}))}\nonumber\\
&&~~~~~~+\frac{D_{h}((u,\omega),(v^{*},z^{*}))}{B((v,z),(v^{*},z^{*}))}
+\frac{D_{h}((v,z),(u^{*},\omega^{*}))}{B((v,z),(v^{*},z^{*}))}.
\end{eqnarray}
 \indent{\bf Proof.}~
 From (\ref{s2.11}), (\ref{s2.18}), (\ref{s4.1}) and
(\ref{s4.2}) we have
\begin{eqnarray*}
&&A_{h}((u,\omega)-(v,z), (u^{*},\omega^{*})-(v^{*},z^{*}))-\lambda
B((u,\omega)-(v,z),
(u^{*},\omega^{*})-(v^{*},z^{*}))\nonumber\\
&&~~~=A_{h}((u,\omega),(u^{*},\omega^{*}))+A_{h}((v,z),(v^{*},z^{*}))-A_{h}((u,\omega),(v^{*},z^{*}))\nonumber\\
&&~~~~~~-A_{h}((v,z),(u^{*},\omega^{*}))
-\lambda(B((u,\omega),(u^{*},\omega^{*}))+B((v,z),(v^{*},z^{*}))\nonumber\\
&&~~~~~~-B((u,\omega),(v^{*},z^{*}))
-B((v,z),(u^{*},\omega^{*})))\\
&&~~~=\lambda
B((u,\omega),(u^{*},\omega^{*}))+A_{h}((v,z),(v^{*},z^{*}))-B(\lambda
(u,\omega),(v^{*},z^{*}))\nonumber\\
&&~~~~~~-D_{h}((u,\omega),(v^{*},z^{*}))
-B((v,z),\lambda^{*}(u^{*},\omega^{*}))-D_{h}((v,z),(u^{*},\omega^{*}))
\\
&&~~~~~~-\lambda B((u,\omega),(u^{*},\omega^{*}))-\lambda
B((v,z),(v^{*},z^{*}))\nonumber\\
&&~~~~~~+\lambda B((u,\omega),(v^{*},z^{*}))
+\lambda B((v,z),(u^{*},\omega^{*}))\\
&&~~~=A_{h}((v,z),(v^{*},z^{*}))-\lambda
B((v,z),(v^{*},z^{*}))\nonumber\\
&&~~~~~~+D_{h}((u,\omega),(v^{*},z^{*}))+D_{h}((v,z),(u^{*},\omega^{*})),
\end{eqnarray*}
dividing $B((v,z),(v^{*},z^{*}))$ in both side of the above we
obtain the desired conclusion.~~~$\Box$\\

\indent{\bf Theorem 4.3.}~ Assume that the conditions of Theorem 4.2
are valid, and $R(E), R(E^{*})$ $\subset W^{3,p}(\Omega)$ $(p\in
\mathbb{S})$, the ascent $\alpha$ of $\lambda$ is equal to 1, for
any $(v,z) \in \mathbf{H}_{h}+\mathbf{H}$, (\ref{s3.3}) and
(\ref{s3.4}) hold, then
\begin{eqnarray}\label{s4.19}
|\lambda_{h}-\lambda|\lesssim h^{2+2(\frac{1}{2}-\frac{1}{p})d}.
\end{eqnarray}
\indent{\bf Proof.}~
 From $\alpha=1$, we know $R(E^{*})$ is the
space of eigenfunctions associated with $\lambda^{*}$. Let
$(u,\omega)$ and $(u_{h},\omega_{h})$ satisfy (\ref{s4.10}) and
(\ref{s4.15}), since $(u,\omega)\in R(E)$, 
 $\|(u,\omega)\|_{A}=1$ ,
 Define
\begin{eqnarray*}
f((v,z))=A(E(v,z), (u,\omega)),~~~\forall (v,z)\in \mathbf{H}.
\end{eqnarray*}
Since  for all $(v,z)\in\mathbf{H}$ one has
\begin{eqnarray*}
&&|f((v,z))|=|A(E(v,z), (u,\omega))|\le
\|E(v,z)\|_A \|(u,\omega)\|_{A}\\
&&~~~ \lesssim \sqrt{\lambda}\|E(v,z)\|_{\mathbf{H}_{1}}  \lesssim
\|E\|_{\mathbf{H}_{1}}\|(v,z)\|_A,
\end{eqnarray*}
$f$ is a linear and bounded functional on $\mathbf{H}$ and
$\|f\|_{A}\lesssim \|E\|_{\mathbf{H}_{1}}$.
 Using Riesz Theorem, we know  there exists
$(u^{*},\omega^{*})\in \mathbf{H}$ satisfying $\|(u^{*},
\omega^{*})\|_{A}=\|f\|_{A}$ and
\begin{eqnarray}\label{s4.20}
A((v,z), (u^{*}, \omega^{*}))=A(E(v,z), (u,\omega)).
\end{eqnarray}
For any $ (v,z)\in \mathbf{H}$, notice $E(I-E)(v,z)=0$,
\begin{eqnarray*}
&&A((v,z),(\lambda^{*-1}-T^{*})(u^{*}, \omega^{*})
)=A((\lambda^{-1}-T)^{\alpha}(v,z),(u^{*}, \omega^{*}) )\\
&&~~~=A((\lambda^{-1}-T)E(v,z),(u^{*}, \omega^{*})
)+A((\lambda^{-1}-T)(I-E)(v,z),(u^{*}, \omega^{*}) )=0,
\end{eqnarray*}
i.e., $(\lambda^{*-1}-T^{*})(u^{*}, \omega^{*})=0$, hence $(u^{*},
\omega^{*})\in R(E^{*})$. From (\ref{s2.32})-(\ref{s2.35}) we have
\begin{eqnarray}\label{s4.21}
&&\|(u^{*},\omega^{*})-I_{h}(u^{*},\omega^{*})\|_{h}\lesssim
h^{1+(\frac{1}{2}-\frac{1}{p})d},\\\label{s4.22}
&&\|(u^{*},\omega^{*})-I_{h}(u^{*},\omega^{*})\|_{\mathbf{H_{1}}}\lesssim
h^{2+(\frac{1}{2}-\frac{1}{p})d},\
\end{eqnarray}
By (\ref{s4.20}) we have
\begin{eqnarray}\label{s4.23}
A((u,\omega), (u^{*}, \omega^{*}))=A(E(u,\omega), (u,
\omega))==A(E(u,\omega), (u, \omega))=1.
\end{eqnarray}
Then, from (\ref{s4.15}), (\ref{s4.21}) and (\ref{s4.23}), when $h$
is small enough, $|A_{h}((u_{h},\omega_{h}),
I_{h}(u^{*},\omega^{*}))|$ has a positive lower bound uniformly with
respect to $h$, thus there is a positive constant $C_{0}$
independent of $h$ such that
\begin{eqnarray}\label{s4.24}
|B((u_{h},\omega_{h}),I_{h}(u^{*},\omega^{*}))|=|\lambda_{h}^{-1}A_{h}((u_{h},\omega_{h}),
I_{h}(u^{*},\omega^{*}))|\geq C_{0}.
\end{eqnarray}
In (\ref{s4.18}), let $(v,z)=(u_{h},\omega_{h}),
(v^{*},z^{*})=I_{h}(u^{*},\omega^{*})$, noting that\\
$$\lambda_{h}=A((u_{h},\omega_{h}),I_{h}(u^{*},\omega^{*}))/B((u_{h},\omega_{h}),I_{h}(u^{*},\omega^{*})),$$
then
\begin{eqnarray}\label{s4.25}
&&|\lambda_{h}-\lambda|\lesssim
\|(u,\omega)-(u_{h},\omega_{h})\|_{h}
\|(u^{*},\omega^{*})-I_{h}(u^{*},\omega^{*})\|_{h}\nonumber\\
&&~~~~~~+\|(u,\omega)-(u_{h},\omega_{h})\|_{\mathbf{H}_{1}}
\|(u^{*},\omega^{*})-I_{h}(u^{*},\omega^{*})\|_{\mathbf{H}_{1}}\nonumber\\
&&~~~~~~+|D_{h}((u,\omega),I_{h}(u^{*},\omega^{*}))|
+|D_{h}((u_{h},\omega_{h}),(u^{*},\omega^{*}))|.
\end{eqnarray}
From (\ref{s3.3}) and (\ref{s4.21}),
\begin{eqnarray}\label{s4.26}
&&|D_{h}((u,\omega),I_{h}(u^{*},\omega^{*}))|=|D_{h}((u,\omega),I_{h}(u^{*},\omega^{*})-(u^{*},\omega^{*}))|\nonumber\\
&&~~~\lesssim
h^{1+(\frac{1}{2}-\frac{1}{p})d}\|(u,\omega)\|_{3,p}\|I_{h}(u^{*},\omega^{*})-(u^{*},\omega^{*})\|_{h}\nonumber\\
&&~~~\lesssim
h^{2+2(\frac{1}{2}-\frac{1}{p})d}\|(u,\omega)\|_{3,p}\|(u^{*},\omega^{*})\|_{3,p},
\end{eqnarray}
and from (\ref{s3.4}) and (\ref{s4.15}),
\begin{eqnarray}\label{s4.27}
&&|D_{h}((u_{h},\omega_{h}),(u^{*},\omega^{*}))|=|D_{h}((u_{h},\omega_{h})-(u,\omega),(u^{*},\omega^{*}))|\nonumber\\
&&~~~\lesssim
h^{2+2(\frac{1}{2}-\frac{1}{p})d}\|(u,\omega)\|_{3,p}\|(u^{*},\omega^{*})\|_{3,p}.
\end{eqnarray}
Substituting (\ref{s4.15}), (\ref{s4.21}), (\ref{s4.22}),
(\ref{s4.26}) and (\ref{s4.27}) into (\ref{s4.25}), we get
(\ref{s4.19}).~~~$\Box$\\

\indent {\bf Remark 4.1.}~~Using the same argument as in this
section
 we can prove the error estimates of finite element approximation for the dual problem
 (\ref{s2.18}): Let $R(E^{*})\subset W^{3,p}(\Omega)$ $(p\in
\mathbb{S})$, then
\begin{eqnarray}\label{s4.28}
 &&\|(u_{h}^{*},\omega_{h}^{*})-(u^{*},\omega^{*})\|_{h}\lesssim
h^{\frac{1}{\alpha}+(\frac{1}{2}-\frac{1}{p})\frac{d}{\alpha}},\\\label{s4.29}
&&\|(u_{h}^{*},\omega_{h}^{*})-(u^{*},\omega^{*})\|_{\mathbf{H}_{1}}\lesssim
h^{\frac{1}{\alpha}+(\frac{1}{2}-\frac{1}{p})\frac{d}{\alpha}}.
\end{eqnarray}


\indent{\bf Theorem 4.4.}~ Assume that the conditions of Theorem 4.2
are valid, and $D\subset \mathbb{R}^{2}$, $R(E)\subset
H^{4}(\Omega)$, $n$ is a constant. Let
  $S^{h}$ be the Adini element space defined on the uniform rectangle mesh. Then
\begin{eqnarray}\label{s4.30}
&&\|(u_{h},\omega_{h})-(u,\omega)\|_{h}\lesssim h^{\frac{2}{\alpha}},\\
\label{s4.31} &&|\lambda_{h}-\lambda|\lesssim
h^{\frac{2}{\alpha}},\\\label{s4.32}
&&|(\frac{1}{q}\sum\limits_{i=1}^{q}\lambda_{i,h}^{-1})^{-1}-\lambda|\lesssim
h^{2}.
\end{eqnarray}
\indent{\bf Proof.}~ From line 11 on page 304 in \cite{ciarlet}, we
know for Adini element, the third term on the right-hand of
(\ref{s3.9}) is equal to $0$, thus we obtain
\begin{eqnarray}\nonumber
&&D_{h}((\psi,\varphi),(v,z))
=-\sum\limits_{\kappa\in\pi_{h}}\int\limits_{\partial\kappa}
\frac{1}{n-1} \Delta \psi
 \nabla \overline{v}\cdot\gamma ds\\\label{s4.33}
&&~~~+\mu_{1}\sum\limits_{\kappa\in\pi_{h}}\int\limits_{\partial\kappa}
\{\frac{\partial^{2}\psi}{\partial x_{1}^{2}}\frac{\partial
\overline{v}}{\partial
x_{2}}\gamma_{2}+\frac{\partial^{2}\psi}{\partial
x_{2}^{2}}\frac{\partial \overline{v}}{\partial x_{1}}\gamma_{1}\}ds
+0.
\end{eqnarray}
Comparing it with the consistency term of the clamped plate bending
problem (see (50.7) of \cite{ciarlet}), from \cite{lascaux} we can
deduce
\begin{eqnarray}\label{s4.34}
&&|D_{h}((\psi,\varphi),(v,z))| \lesssim
h^{2}\|\psi\|_{4}\|v\|_{h},~~~\forall (v,z)\in
\mathbf{H}_{h},\\\label{s4.35}
&&\|(\psi,\varphi)-(\psi_{h},\varphi_{h})\|_{h}\lesssim
h^{2}(\|\psi\|_{4}+\|\varphi\|_{2}).
\end{eqnarray}
Thus we have
\begin{eqnarray}\label{s4.36}
\|(T-T_{h})|_{R(E)}\|_{\mathbf{H}_{1}} =\sup\limits_{(u,\omega)\in
R(E),\|(u,\omega)\|_{\mathbf{H}_{1}}=1}\|(T-T_{h})(u,\omega)\|_{\mathbf{H}_{1}}\nonumber\\
\lesssim \sup\limits_{(u,\omega)\in
R(E),\|(u,\omega)\|_{\mathbf{H}_{1}}=1}h^{2}\|T(u,\omega)\|_{H^{4}(\Omega)\times
H^{2}(\Omega)} \lesssim h^{2}.
\end{eqnarray}
Substituting (\ref{s4.36}) into (\ref{s4.11}) and (\ref{s4.12}) we
get (\ref{s4.31}) and (\ref{s4.32}), respectively.\\
\indent By the way to show (\ref{s4.15}), we can prove
(\ref{s4.30}).~~~$\Box$\\

\indent The literature \cite{yang1} proved that the order of
convergence is just 2 for the Adini finite element eigenvalues for
the clamped plate vibration problem. Based on  \cite{yang1}, we can
prove the estimate (\ref{s4.31}) and (\ref{s4.32}) are optimal and
cannot be improved further.

\section{Numerical Experiment}
\indent In this section, we will report  some numerical experiments
for non-conforming finite element discretizations
to validate our theoretical results.\\
\indent We use Matlab 2012a to solve (\ref{s2.1})-(\ref{s2.4}) on a
Lenovo G480 PC with 4G memory. Our program is compiled under the
package of iFEM \cite{chen}. \\
\indent Let $\{\xi_j\}_{j=1}^{N_h}$ be a basis of ${S}^h$ and
$u_h=\sum_{j=1}^{N_h}u_j\xi_j,\omega_h=\sum_{j=1}^{N_h}\omega_j\xi_j$.
Denote $\overrightarrow{u}=(u_1,\cdots,u_{N_h})^T$ and
$\overrightarrow{\omega}=(\omega_1,\cdots,\omega_{N_h})^T$. To
describe our algorithm, we specify the following $N_h\times N_h$ matrices in  the
discrete case.
\begin{center} \footnotesize
\begin{tabular}{lllll}\hline
Matrix&Definition\\\hline $A_h$&$a_{lj}=\sum\int_\kappa
\big\{(\frac{1}{n-1}-\mu_1)\Delta\xi_j\Delta\xi_l +
\mu_1(\frac{\partial^{2}\xi_j}{\partial
x_{1}^{2}}\frac{\partial^{2}\xi_l}{\partial x_{1}^{2}}
+2\frac{\partial^{2}\xi_j}{\partial x_{1}\partial
x_{2}}\frac{\partial^{2}\xi_l}{\partial x_{1}\partial x_{2}}
+\frac{\partial^{2}\xi_j}{\partial
x_{2}^{2}}\frac{\partial^{2}\xi_l}{\partial x_{2}^{2}})\big\}$\\
$B_h$&$b_{lj}=\int_D \big\{\nabla(\frac{1}{n-1}\xi_j)\cdot\nabla
\xi_l +
 \nabla \xi_j\cdot \nabla(\frac{n}{n-1}\xi_l) \big\}dx$\\
$C_h$&$c_{lj}=-\int_D \frac{n}{n-1} \xi_j \xi_ldx$\\
$G_h$&$g_{lj}=\int_D  \xi_j \xi_ldx$\\
\hline
\end{tabular}
\end{center}
\noindent where $N_h=\mathrm{dim}(S^h)$. Then (\ref{s2.23}) can be
written as a generalized eigenvalue problem
\begin{eqnarray} \label{s5.1}
\left(
\begin{array}{lcr}
A_h&0\\
0&G_h
\end{array}
\right) \left (
\begin{array}{lcr}
\overrightarrow{u}\\
\overrightarrow{\omega}
\end{array}
\right)=\lambda_h  \left(
\begin{array}{lcr}
B_h&C_h\\
G_h&0
\end{array}
\right) \left (
\begin{array}{lcr}
\overrightarrow{u}\\
\overrightarrow{\omega}
\end{array}
\right).
\end{eqnarray}
\indent Note that in (\ref{s5.1}) $A_{h}$ is a positive definite
Hermitian matrix, and $G_{h}$ can be equivalently replaced by the
identity matrix $I_h$, which will lead to two sparser coefficient
matrices with good structure.  Based on this fact, 
we use the sparse matrix eigenvalue solver $eigs$ to compute the
numerical eigenvalues and
the resulting numerical eigenvalues are  ideal. \\
\indent We consider the model problems (\ref{s2.1})-(\ref{s2.4})
with the refraction index $n=8+x_1-x_2$ and $n=16$ on the unit
square $(0,1)^2$, L-shaped $(-1,1)^2\backslash\big([0,1)\times(-1,0]\big)$,
triangle whose vertices are given by
$(-\frac{\sqrt3}{2},-\frac{1}{2})$,
$(\frac{\sqrt3}{2},-\frac{1}{2})$ and
 $(0,1)$, and disk with radius $\frac{1}{2}$ and center $(0,0)$. We adopt the Morley-Zienkiewicz(MZ)
  element and Adini element to compute the transmission eigenvalues on quasi-uniform meshes.
   The Morley-Zienkiewicz element was put  forward in \cite{shi}, and its finite element space is defined
   as:\\
$S^h=\{v\in V_h|v~ {and}~\nabla v ~ {vanish~ at~ all~ vertices~ on}~ \partial \Omega, ~ and ~over~ any ~face ~F\\
   \subset\partial \Omega, ~ { the~ mean~ value ~of~} \frac{\partial v}{\partial\nu_F} ~ {vanishes}\}$,\\
 where\\
  $V_h=\{v\in L^2(\Omega):
 v|_\kappa\in P_\kappa,~\forall \kappa\in \pi_h;~v~ {and}~\nabla  v ~are ~continuous~ at ~all ~vertices\\
  ~of  ~\pi_h,   {and~over~each~interelement~face~F~of}~\pi_h,~the~jump ~of~the ~mean ~value\\
   ~of~ \frac{\partial v}{\partial\nu_F} ~is ~zero\}$,\\
   $P_\kappa=P_3''(\kappa)+span\{l_1^{2}l_2\cdots l_{d+1},l_1 l_2^2\cdots l_{d+1},\cdots,l_1 l_2\cdots l_{d+1}^2\}$
  with $P_3''(\kappa)$ being the Zienkiewicz element shape function
  space and $l_i(i=1,\cdots,d+1)$ being the barycentric coordinates.\\
\indent In our computation, for all the domains mentioned above we
set $\mu_1=\frac{1}{9}$ when the refraction index $n=8+x_1-x_2$ and
$\mu_1=\frac{1}{15}$ when the refraction index $n=16$.
 The associated numerical eigenvalues computed by the MZ element and Adini element
 are listed partially in Tables 1-2 and Table 3, while the error curves of these numerical eigenvalues
 whose slopes are computed by  procedure of curve fitting are depicted in Figures
 1-3.  \\
\indent For reading conveniently,  in our tables and figures we use
the  notation  $k^{\Omega}_{j,h}=\sqrt{\lambda^{\Omega}_{j,h}}$ to
denote the $jth$ eigenvalue on the domain $\Omega=S, L, T, D$
obtained by (\ref{s2.23}) on $\pi_{h}$, where the symbols $S, L, T,
D$ denote the domains  square, L-shaped, triangle and disk,
respectively.
\\
\indent  It is seen from Figures 1-3 that the convergence orders of
the numerical eigenvalues on the unit square, triangle and disk
computed by the two elements are around 2, which coincides with the
theoretical result. Nevertheless, the convergence orders on the
L-shaped domain of the numerical eigenvalues
$k_{1,h},k_{2,h},k_{5,h},k_{6,h}$ with $n=8+x_1-x_2$ and
$k_{1,h},k_{3,h}$ with $n=16$ are less than 2 (see Figures 1-2).
This fact suggests that the eigenfunctions corresponding to these
eigenvalues  on the L-shaped domain do have singularities to
different degrees.
Numerical results  indicates our discretizations by the MZ element and the Adini
element are efficient and consistent with theoretical analysis.

\begin{table}
\caption{The eigenvalues obtained by MZ element, $n=8+x_1-x_2$.}
\begin{center} \footnotesize
\begin{tabular}{ccccccc}\hline
$j$&$h$ & $k^S_{j,h}$
&$k^L_{j,h}$&$k^T_{j,h}$&$h$&$k^D_{j,h}$\\\hline
1&  $\frac{\sqrt2}{32}$&        2.8218574&  2.3035843&  2.7388174&0.025&2.9775769\\
1&  $\frac{\sqrt2}{64}$&      2.8220628&   2.3028188&   2.7389418&0.012&2.9771919\\
1&  $\frac{\sqrt2}{128}$&    2.8221545&  2.3024576&    2.7389765&0.006&2.9771000\\
2&  $\frac{\sqrt2}{32}$&  3.5381161&  2.3953577&  3.2915472&0.025&3.7774560\\
2&  $\frac{\sqrt2}{64}$&   3.5384282& 2.3955964&   3.2917188&0.012&3.7770363\\
2&  $\frac{\sqrt2}{128}$& 3.5386203&  2.3956673&   3.2917696&0.006&3.7769414\\
5,6&  $\frac{\sqrt2}{32}$&    4.4959659 &    2.9255876&     4.1666454&0.025&  4.8741035\\
&&$\pm$0.8714721i&$\pm$0.5654338i&$\pm$0.7836432i&&$\pm$0.8760355i\\
5,6&  $\frac{\sqrt2}{64}$&     4.4963441&   2.9248145&     4.1666973&0.012&  4.8733986\\
&&$\pm$0.8714728i&$\pm$0.5650876i&$\pm$0.7836699i&&$\pm$0.8758772i\\
5,6&  $\frac{\sqrt2}{128}$&   4.4964963&    2.9244878&     4.1667103&0.006&  4.8732345 \\
&&$\pm$0.8714802i&$\pm$0.5648487i&$\pm$0.7836780i&&$\pm$0.8758363i\\
 \hline
\end{tabular}
\end{center}
\end{table}

\begin{table}
\caption{The eigenvalues obtained by MZ element, $n=16$.}
\begin{center} \footnotesize
\begin{tabular}{ccccc|ccc}\hline
$j$&$h$ & $k^S_{j,h}$
&$k^L_{j,h}$&$k^T_{j,h}$&$j$&$h$&$k^D_{j,h}$\\\hline
1&  $\frac{\sqrt2}{32}$&     1.8795675&   1.4775023&    1.8184414&1&0.025&1.9883914\\
1&  $\frac{\sqrt2}{64}$&     1.8795717&   1.4767526&   1.8184573&1&0.012&1.9880919\\
1&  $\frac{\sqrt2}{128}$&    1.8795854&  1.4764066&   1.8184622&1&0.006&1.9880191\\
2&  $\frac{\sqrt2}{32}$& 2.4440863&  1.5696996&  2.2870296&2,3&0.025&2.6134315\\
2&  $\frac{\sqrt2}{64}$& 2.4441734& 1.5697172&   2.2870557&2,3&0.012&2.6130503\\
2&  $\frac{\sqrt2}{128}$& 2.4442186&  1.5697237&  2.2870651&2,3&0.006&2.6129596\\
3&  $\frac{\sqrt2}{32}$&  2.4442285&  1.7053198&   2.2870296&13,14&0.049& 4.9056584\\
3&  $\frac{\sqrt2}{64}$&  2.4441893&  1.7051917&   2.2870557&&& $\pm$0.5787253i\\
3&  $\frac{\sqrt2}{128}$& 2.4442212&  1.7051196&   2.2870651&13,14&0.025& 4.9018623 \\
4&  $\frac{\sqrt2}{32}$&  2.8667518&  1.7830953&  2.8375736&&&$\pm$0.5781361i\\
4&  $\frac{\sqrt2}{64}$&  2.8664156&  1.7831002&  2.8376056&13,14& 0.006&  4.9009219\\
4&  $\frac{\sqrt2}{128}$& 2.8664256&  1.7831114&  2.8376222&&&
$\pm$0.5781031i
 \\\hline
\end{tabular}
\end{center}
\end{table}

\begin{table}
\caption{The eigenvalues obtained by Adini element on the unit
square.}
\begin{center} \footnotesize
\begin{tabular}{cccccc}\hline
$h$&$j$ & $k^S_{j,h}$($n=8+x_1-x_2$)&
$j$&$k^S_{j,h}$($n=16$)\\\hline
  $\frac{\sqrt2}{32}$&1&     2.8178682& 1&  1.8778418\\
  $\frac{\sqrt2}{64}$&1&     2.8211011& 1&  1.8791512\\
  $\frac{\sqrt2}{128}$&1&    2.8219168& 1&  1.8794810\\
  $\frac{\sqrt2}{32}$&2& 3.532859351&     2,3& 2.4413924\\
  $\frac{\sqrt2}{64}$&2& 3.537222143&     2,3& 2.4435179\\
  $\frac{\sqrt2}{128}$&2& 3.538327097&    2,3& 2.4440561\\
  $\frac{\sqrt2}{32}$&5,6& 4.4949831$\pm$0.8710067i&     4&   2.8588866\\
  $\frac{\sqrt2}{64}$&5,6& 4.4961529$\pm$0.8713583i&     4&   2.8645286\\
  $\frac{\sqrt2}{128}$&5,6& 4.4964517$\pm$0.8714506i&    4&   2.8659601\\
\hline
\end{tabular}
\end{center}
\end{table}

 \begin{figure}
\includegraphics[width=0.4\textwidth]{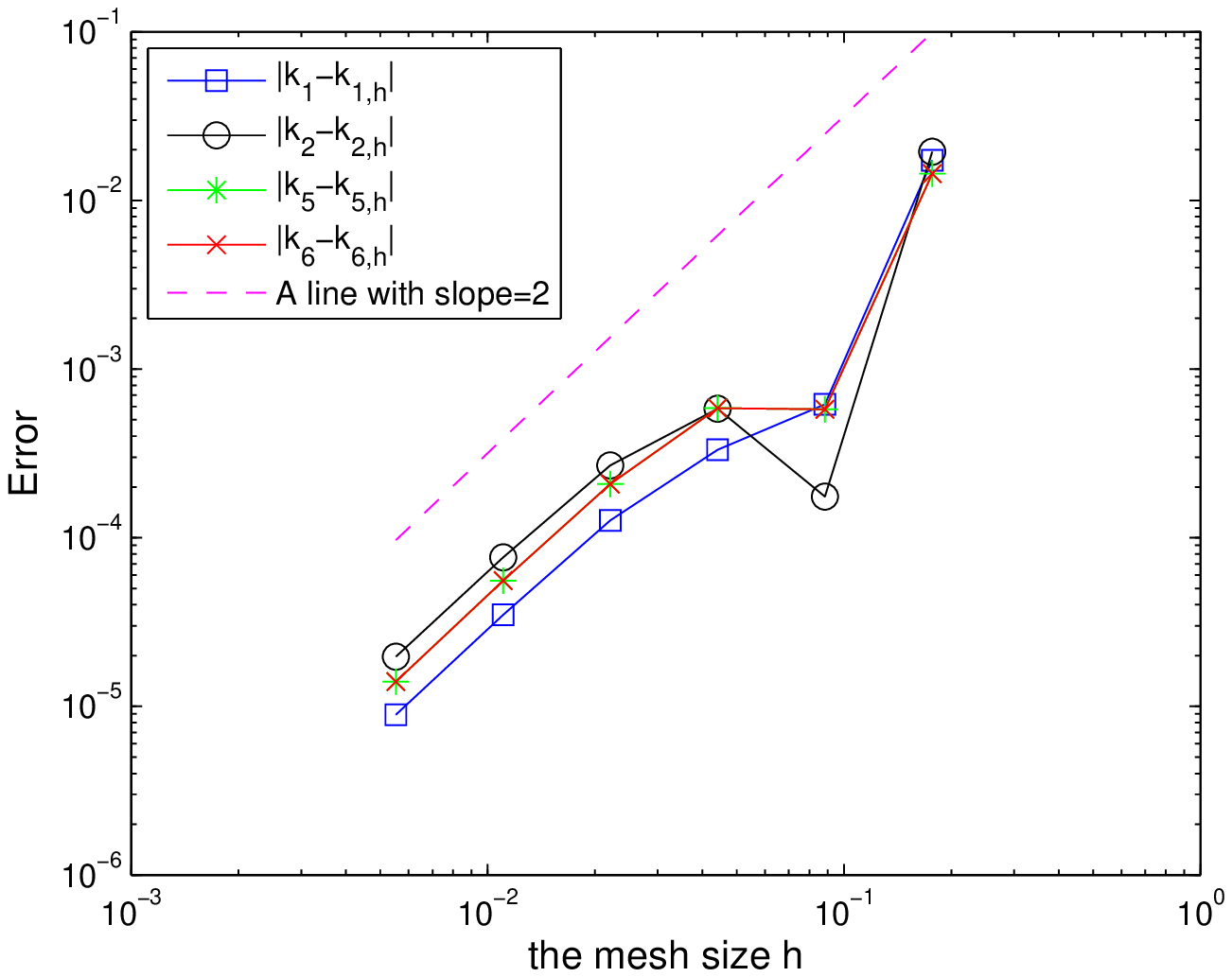}
 \includegraphics[width=0.4\textwidth]{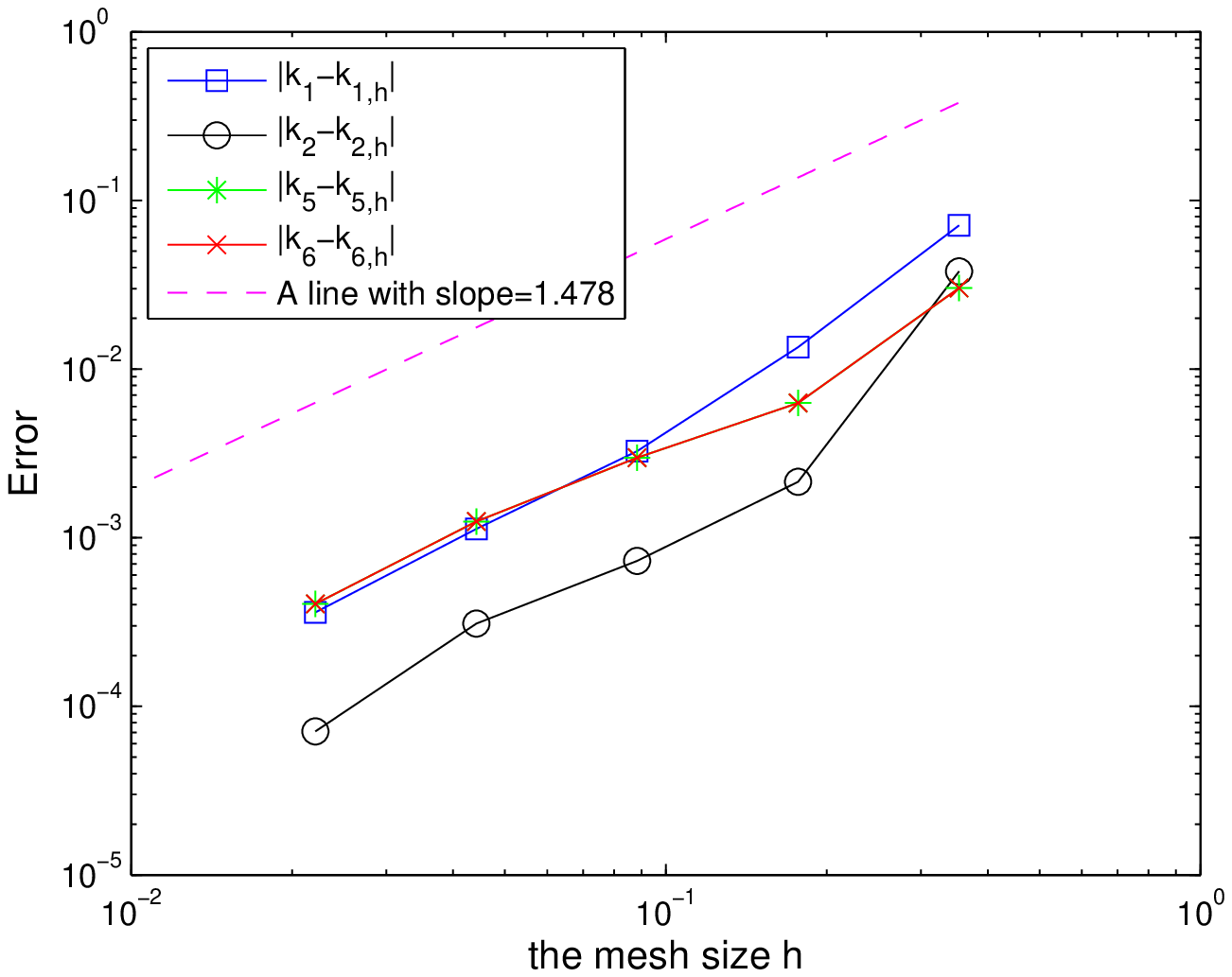}\\
 \includegraphics[width=0.4\textwidth]{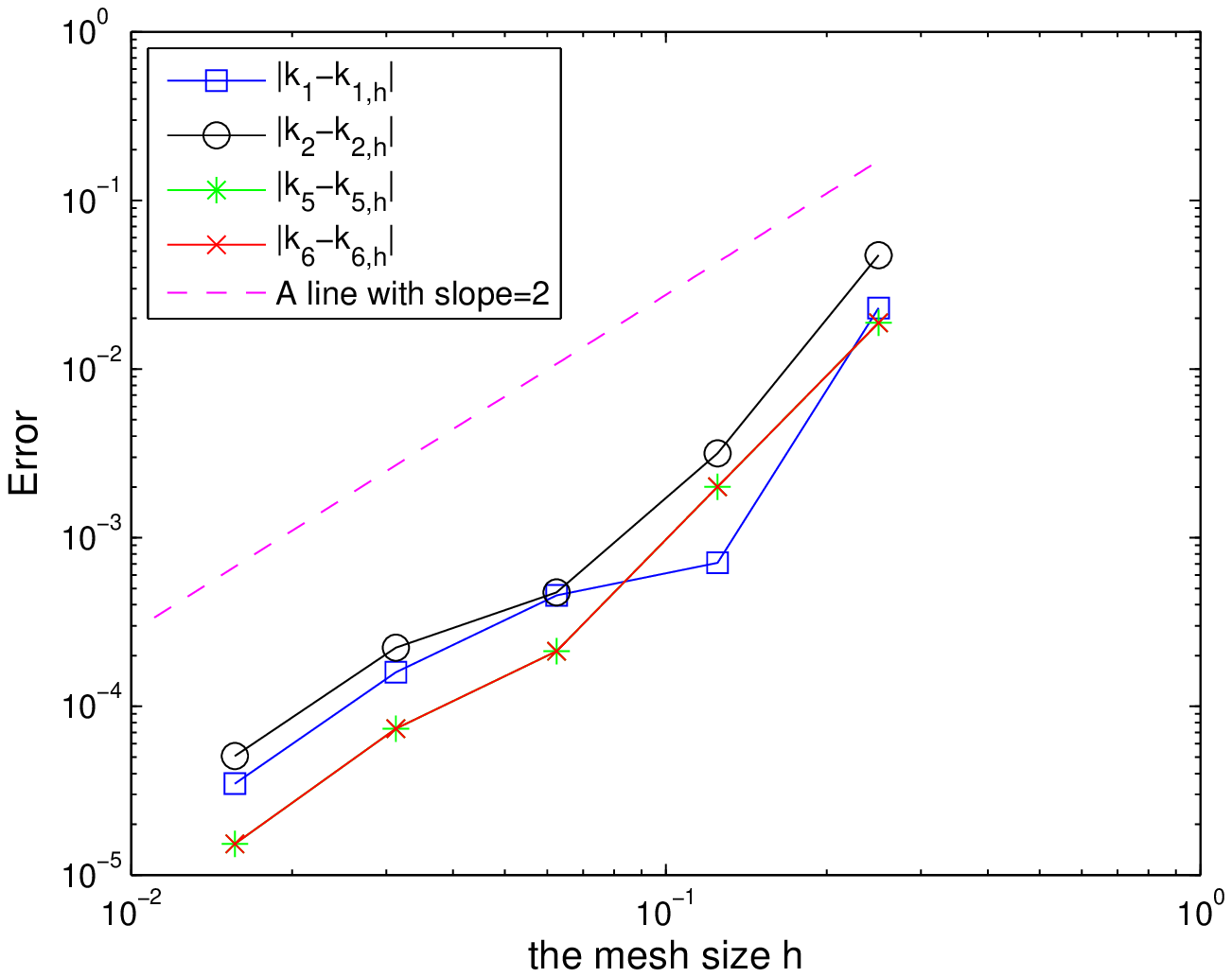}
  \includegraphics[width=0.4\textwidth]{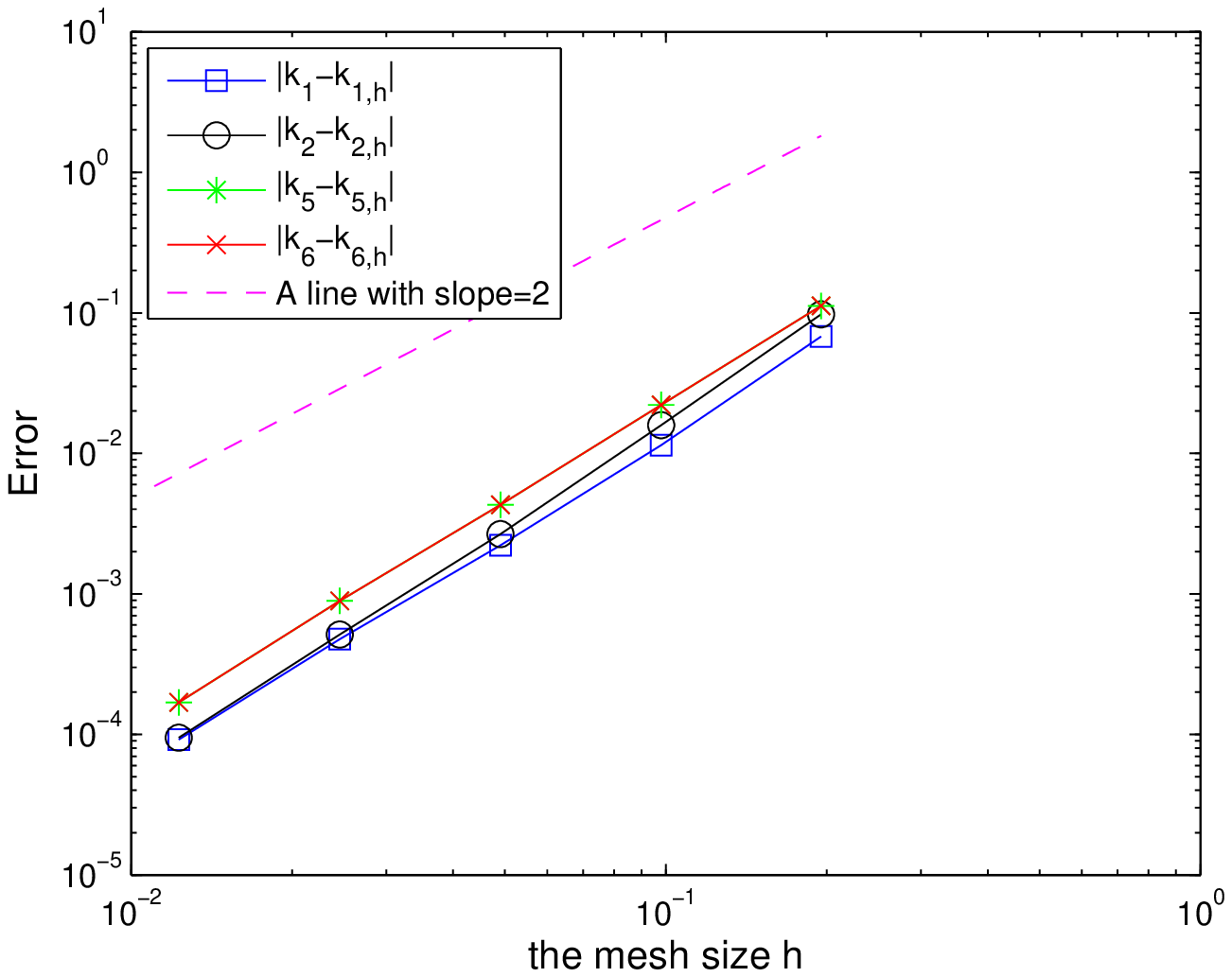}
\caption{{ Error curves  computed by MZ element with $n=8+x_1-x_2$
on the unit square (left top), on the L-shaped  (right top), on the
triangle  (left bottom), on the disk  (right bottom).}}
 \end{figure}

  \begin{figure}
\includegraphics[width=0.4\textwidth]{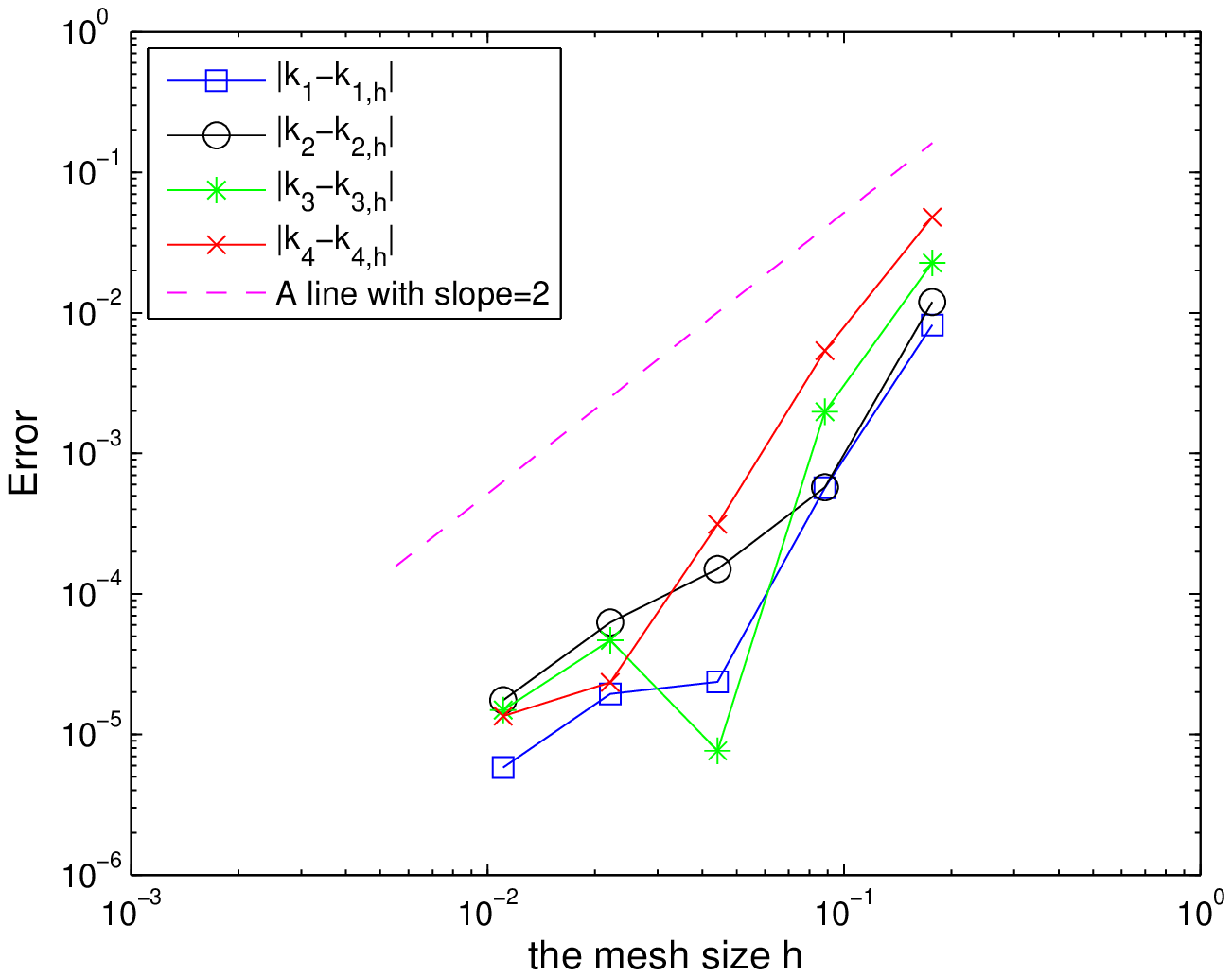}
 \includegraphics[width=0.4\textwidth]{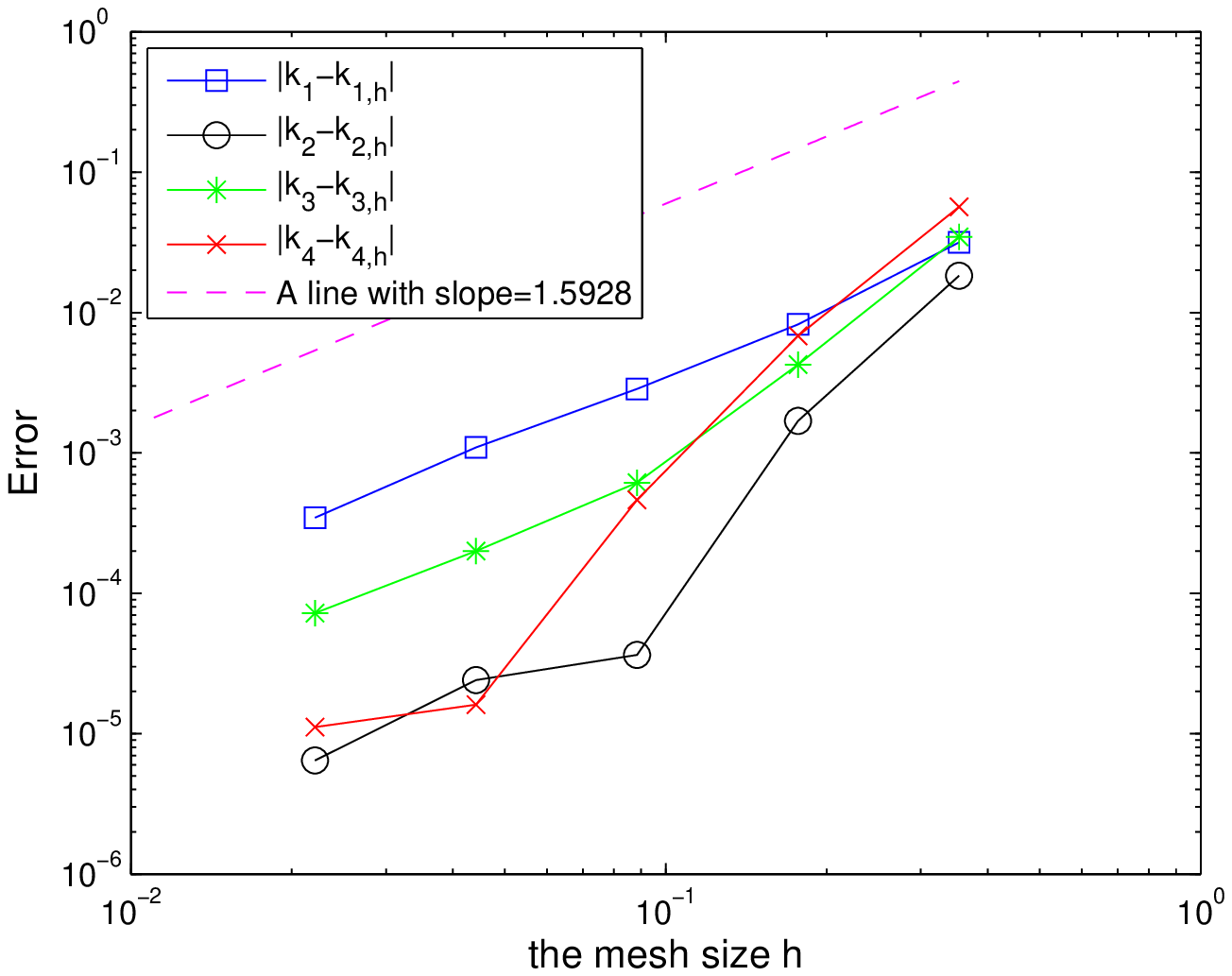}\\
 \includegraphics[width=0.4\textwidth]{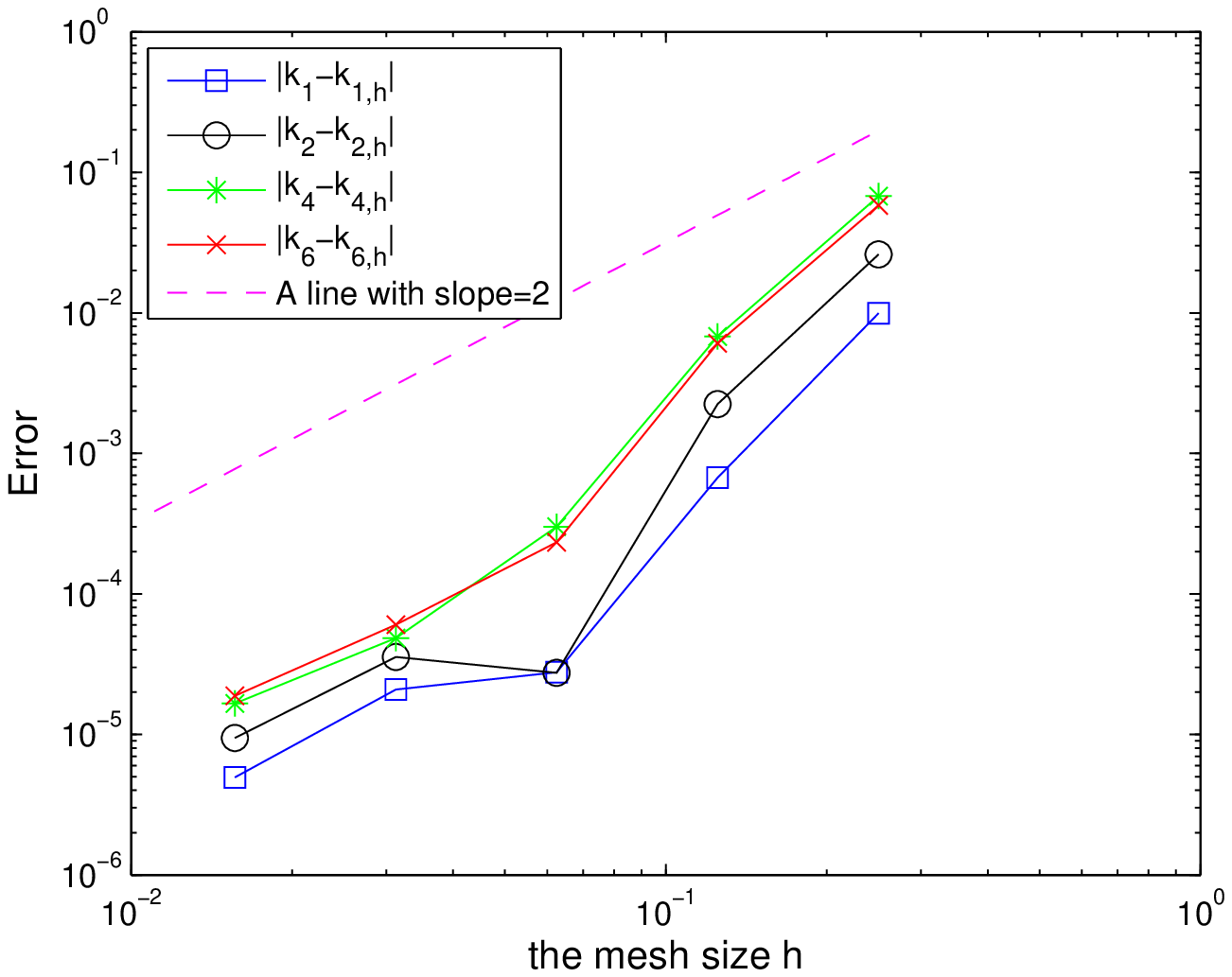}
  \includegraphics[width=0.4\textwidth]{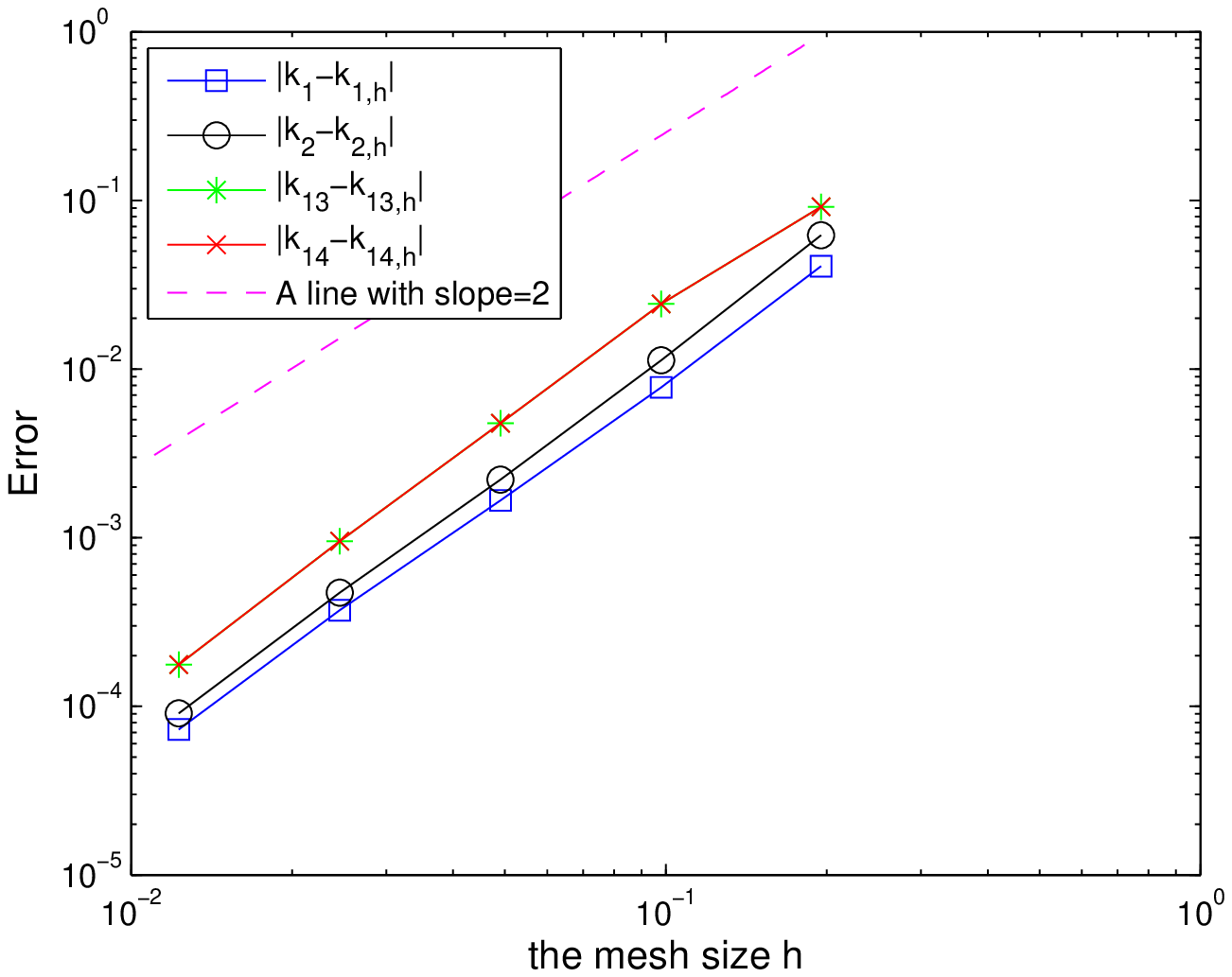}
\caption{{ Error curves  computed by MZ element with $n=16$ on the
unit square (left top), on the L-shaped  (right top), on the
triangle  (left bottom), on the disk  (right bottom).}}
 \end{figure}

   \begin{figure}
\includegraphics[width=0.4\textwidth]{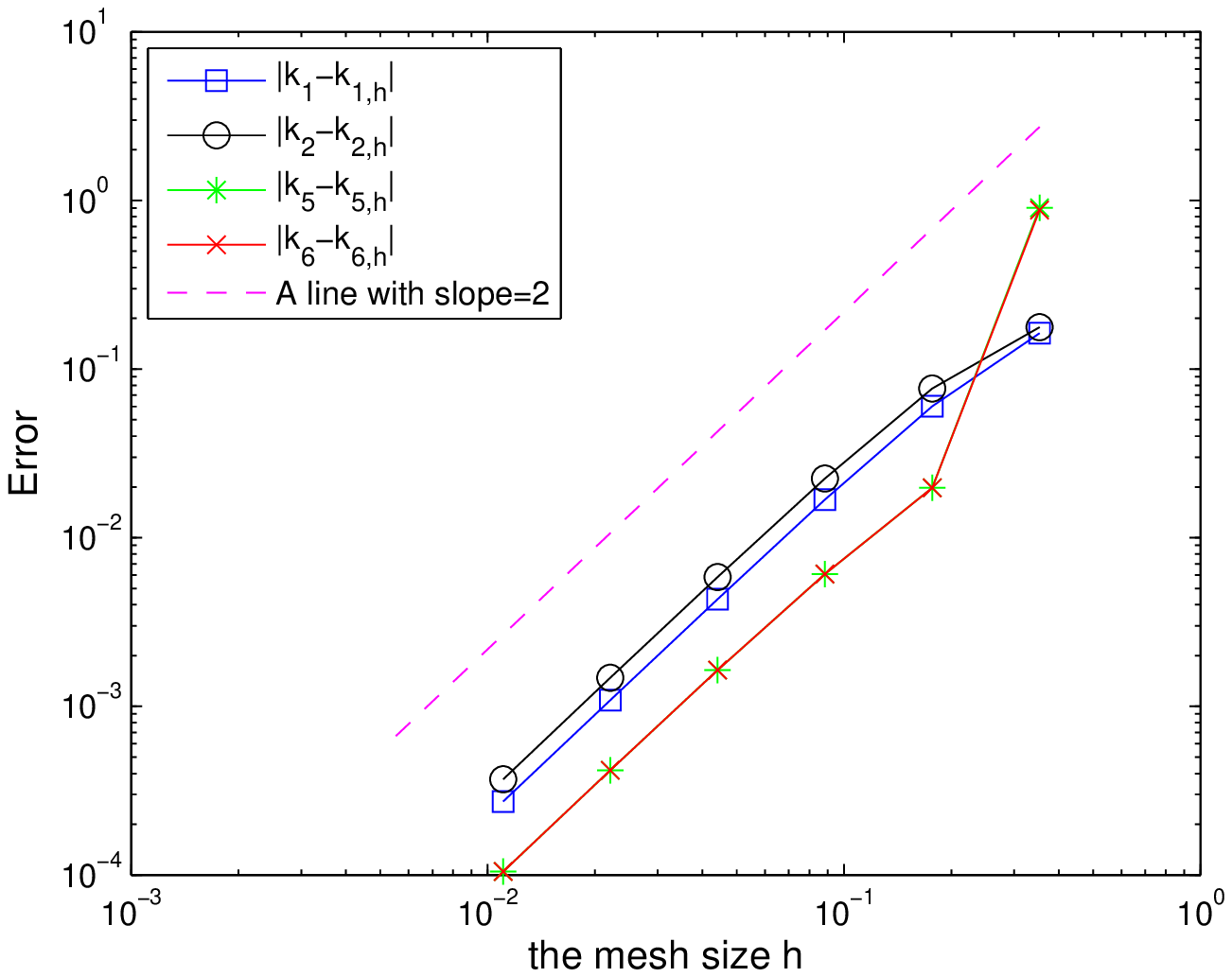}
 \includegraphics[width=0.4\textwidth]{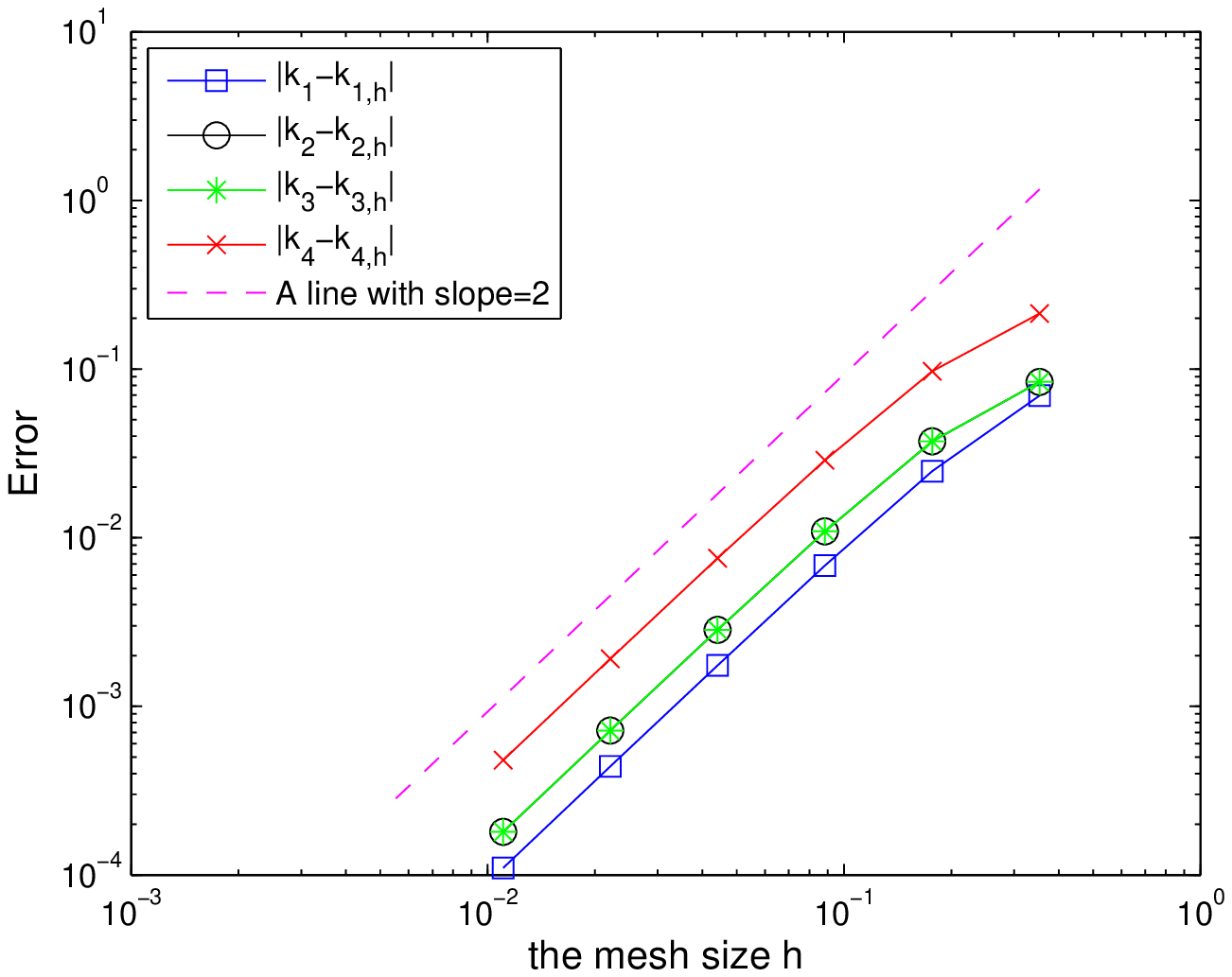}
\caption{{ Error curves  computed by Adini element on the unit
square with $n=8+x_1-x_2$ (left) and with $n=16$ (right).}}
 \end{figure}


\bibliographystyle{amsplain}

\end{document}